\magnification=\magstep{1}
\hoffset0.4 true in
\voffset-30pt
\input amstex
\documentstyle{amsppt}
\nologo
\loadbold
\topmatter
\title {On  boundary properties of solutions of complex vector fields}
\endtitle
\author {S. Berhanu and J.Hounie}
\endauthor
\address
Department of Mathematics, Temple University,
Philadelphia, PA 19122-6094, USA
\endaddress
\email
berhanu\@math.temple.edu
\endemail
\address
Departamento de Matem\'atica, UFSCar,
13.565-905,  S\~ao Carlos, SP,  BRASIL
\endaddress
\email
hounie\@ufscar.dm.br
\endemail
\thanks
 Work supported in part by CNPq, FINEP and FAPESP.
\endthanks

\def\ce{\Bbb C}
\def\erre{\Bbb R}
\def\ccinf{C^\infty_{c}}

\define\<{\langle}

\define\>{\rangle}

\redefine \C{\Cal{C}}
\redefine \D{\Cal{D}}
\redefine \P{\Cal{P}}

\redefine \L{\Cal{L}}

\redefine \|{\Vert}

\abstract{ This work presents results on the boundary properties of solutions
of a complex, planar, smooth vector field $L$. Classical results in the $H^p$
theory of holomorphic functions of one variable are extended to the solutions
of a class of nonelliptic complex vector fields.}
\endabstract
\endtopmatter
\NoBlackBoxes
 \document
 \heading
 {Introduction}
 \endheading

Suppose $h(z)$ is a holomorphic function of one variable defined on some
rectangle
$$Q=(-a,a)\times (0,b)$$
with a weak boundary value at $y=0$. It is well known that if the boundary
value $bh\in L^p(-a,a)$ for some $1\le p\le\infty$ then
\roster
\item for any $0<c<a$, the norms of the traces
$h(.,y)$ in $L^p[-c,c]$  are uniformly bounded as $y\mapsto 0^+$:
$$
\int_{-c}^c |h(x+iy)|^p\,dx\le C,\quad y\searrow0;
$$
\item $h(x+iy)$ converges pointwise and nontangentially to $bh(x)$ for almost every
$x\in(-a,a)$; 
\item $h(x+iy)$ vanishes identically  if $bh$ vanishes on a set of positive measure;
\item conversely, if (1) holds, $bh\in L^p(-c,c)$ for any $0<c<a$.
\endroster
These are just the local versions of very classical properties for
holomorphic functions on the unit disc $\Delta$. Fatou proved in 1906
[F] that any bounded holomorphic function $f$ on $\Delta$ has an a.e.
nontangential limit that cannot vanish identically on an arc of
$\partial\Delta$ unless $f$ is identically zero and that the Poisson
integral of a finite measure has a.e.  nontangential limit. Then Hardy
[Ha] initiated the theory of the spaces $H^p(\Delta)$ in 1915, proving
that the logarithm of the $L^p[-\pi,\pi]$ norm of $\theta\mapsto
f(re^{i\theta})$ is a convex function of $\ln r$, $0<r<1$. The weak
compactness of the unit ball of $L^p$ implies easily the validity of
(4) for $p>1$ but for $p=1$ ---where this argument only yields that
$bf$ is a measure--- it is a consequence of the famous F. and M. Riesz
theorem presented in [RR] in 1916 where it is also shown that any $f\in
H^1(\Delta)$ has an a.e. nontangential limit that cannot vanish
identically on a subset of $\partial\Delta$ of positive measure unless
$f$ is identically zero.

Holomorphic functions are solutions of a complex vector field and in this
paper we explore generalizations of these properties for solutions of more general smooth
complex vector fields in the plane. Our main result, Theorem 3.1 extends the
uniform boundedness of the $L^p$ norms (1) to traces of continuous solutions
of any locally solvable, smooth complex vector field in the plane while Theorems 5.1 and 6.1 address, {\it mutatis mutandis}, 
properties analogous to (2) and (3). The analogue of (4) for the relevant
value $p=1$ was the subject of [BH1].

The uniform control of $L^p$ norms (1) is a distinctive feature of Hardy spaces that was preserved in spite of the extraordinary 
development undergone by the theory  along the century. Present in the
original formulation in connection with boundary values of holomorphic and
harmonic functions, this property was not lost in the real variable definition
of Hardy spaces $H^p(\erre^n)$ in terms of maximal functions by Stein and
Weiss, where the Poisson kernel plays a key role. Indeed, the spaces so
defined coincide with the boundary values of solutions of appropriate elliptic
systems defined on $\erre^n\times(0,\infty)$ for which uniform control in
$t\in(0,\infty)$ of the  $L^p(\erre^n)$ norm holds ---we refer the reader to
the books [St] and [GR] on this subject--- and it seems fair to say that
uniform control of norms for solutions of elliptic equations is not a
surprising fact. On the other hand, uniform control of norms for solutions of
equations which are not necessarily elliptic or even far from elliptic, when
integral representation of solutions in terms of boundary values are not
available, seems new.

Our method of proof takes advantage of  a blend of modern and more classical tools. Among the former we should mention techniques from 
microlocal analysis, specifically the FBI transform  in the fashion developed
in [BCT] and [T1], the  Baouendi-Treves approximation formula [BT] and 
results from the $L^2$ theory of operators with Calder\'on-Zygmund kernels
such as the boundedness of the Cauchy integral and its related maximal
operator. For instance, in the proof of Theorem 3.1, the Baouendi-Treves
formula reduces the problem to the study of (1) for a sequence of holomorphic
functions on planar domains whose boundaries are not smooth due to the
presence of cusps. At this stage, the more classical theory of Jordan domains
with rectifiable boundary as described in chapter 10 of [Du] comes at hand.

For a holomorphic function $h$ defined on $Q$ as above, the existence of a
distribution trace at $y=0$, i.e. the existence of a weak distribution limit
for the traces $h(.,y)$ is equivalent to the property that $h$ be of tempered
growth, i.e. that for some integer $N$, 
$$
h(x+iy)=O(y^{-N})
$$ 
uniformly for $x$ in compact subsets of $(-a,a)$. For continuous solutions of a general, smooth 
complex vector field
$$
L=\frac{\partial}{\partial y}+b(x,y)\frac{\partial}{\partial x}
$$ 
this
equivalence is no longer valid. Indeed, the equivalence fails even for the
subclass of locally integrable vector fields (see the examples in section 1).
In section 1, we prove that if  
$$
L=\frac{\partial}{\partial t}+\sum_{j=1}^{n}a_j(x,t)\frac{\partial}{\partial x_j}
$$
is a smooth complex vector field in $U=B(0,a)\times (0,b)$ in ${\Bbb
R}^{n+1}$, $B(0,a)$ a ball in ${\Bbb R}^n$, $f\in L^1_{\text{loc}}(U)$,
$Lf\in L^1(U)$, and the integrals $$
\int_K|f(x,t)|\,dx=O(t^{-N})
$$
for every $K$ compact in $B(0,a)$, then $\lim_{t\mapsto 0^+}f(x,t)$ exists and
defines a distribution trace at $t=0$. Analogues of this trace result under
more stringent conditions on $f$ were proved in [Br] and [BH1]. In section 2 we
study pointwise convergence of solutions to their $L^p$ boundary values for the
class of locally integrable vector fields. We recall that a nowhere vanishing
planar  vector field $L$ is called locally integrable in an open set $\Omega$
if each $p\in \Omega$ is contained in a neighborhood which admits a smooth
function $Z$ with the properties that $LZ=0$ and the differential $dZ\neq 0$. 
Examples of locally integrable vector fields include nonzero real-analytic
vector fields and locally solvable vector fields. Note however that the class
of locally integrable vector fields is much larger and we refer the reader to
the treatise [T1] for more examples. For solutions of locally integrable
vector fields, as a substitute for radial convergence, we describe curves along which a.e. pointwise convergence holds
on the noncharacteristic portion of the boundary. Examples demonstrate that in
general, one can not get larger sets of approach than these curves. However,
when $L$ is a locally solvable vector field, we will show in section $5$ that
the sets of approach for convergence are open sets at the points where $L$
does not behave like a real vector field in the sense made precise in that
section. Finally, in section 6 we prove a uniqueness result analogous to the
Riesz uniqueness theorem.

\smallskip
\heading {1. A Theorem on the existence of traces}
\endheading

In this section we present conditions that guarantee the existence of a
boundary value for solutions of a complex vector field. It is well known (see
[Ho, Thm.3.1.14]) that if $h$ is holomorphic in a rectangle $Q=(-a,a)\times (0,b)$, then
the traces $h(.,y)$ converge as $y\mapsto 0$ to a distribution $bh(x)$ iff
there exists an integer $N$ such that 
$$|h(x+iy)|=O(y^{-N})$$
uniformly for $x$ in compact sets. In the work [Br] the author generalized one
direction of this result to a smooth complex vector field 
$$
L=\frac{\partial}{\partial t}+\sum_{j=1}^{n}a_j(x,t)\frac{\partial}{\partial x_j},
$$
as follows: 

\proclaim{Theorem} (Theorem {\sl 3.4} in [Br]) Let $X\subseteq {\Bbb
R}^n $ be open, $U$ an open neighborhood  of $X\times\{0\}$ in
${\Bbb R}^{n+1}$, $U_+=U\cap{\Bbb R}_+^{n+1}$.  Let
$L=\frac{\partial}{\partial
t}+\sum_{j=1}^{n}a_j(x,t)\frac{\partial}{\partial x_j}$,
$a(x,t)\in C^{\infty}$ on $X\cup U_+$. Let $f\in C^1(U_+)$ such that
\roster
 \item "i)" $Lf\in L^{\infty}(U_+)$;
\item"ii)"  for any  compact set $K\subset\subset X$ there exists $N=N(K)\in {\Bbb N}$, and $C=C(K)>0$ such that
$$
  |f(x,t)| \le \frac{C}{t^N} \quad, \text{ and } |D_xf(x,t)| \le \frac{C}{t^N}
 $$
\endroster
 Then $\lim_{t\to 0}f(x,t)=bf$
exists in ${\D}'(X)$.
\endproclaim
This result in [Br] was improved in our work [BH1] by dropping the growth
condition on $D_xf(x,t)$ and weakening the regularity of $f$ to continuity.
In both references, the function $f$ was assumed to be of tempered growth as
$t\mapsto 0^+$. In the next theorem, we relax this condition and assume
instead that the integrals of $|f(.,t)|$ over compact subsets are of
tempered growth. We also weaken the regularity assumptions on $f$ and
$Lf$. This stronger trace result allows us to improve the F. and M. Riesz
theorem we proved in [BH1] (see Corollary 1.3).

\proclaim{Theorem 1.1}  Let $X\subseteq {\Bbb
R}^n $ be open, $U$ an open neighborhood  of $X\times\{0\}$ in
${\Bbb R}^{n+1}$, $U_+=U\cap{\Bbb R}_+^{n+1}$.  Let
$L=\frac{\partial}{\partial
t}+\sum_{j=1}^{n}a_j(x,t)\frac{\partial}{\partial x_j}$,
$a(x,t)\in C^{\infty}$ on $X\cup U_+$. Let $f$ be a locally integrable function on $U_+$ such that
 \itemitem{i)} $Lf\in L^1(U_+)$;
\itemitem{ii)}  for any  compact set $K\subset\subset X$ there exists $N=N(K)\in {\Bbb N}$, and $C=C(K)>0$ such that
$$
  \int_K |f(x,t)|\,dx \le \frac{C}{t^N} \quad, \text{ as } t\to 0.
 $$

\noindent
 Then $\lim_{t\to 0}f(x,t)=bf$
exists in ${\D}'(X)$. Furthermore, if $X\times (0,T]\subseteq
U_+$, then the distributions $\{f(.,t):0\leq t \leq T\}$ are
uniformly bounded in ${\D}'(X)$.
\endproclaim

Before we prove Theorem 1.1, we present two examples where this theorem
can be applied. In both cases, the solution $f$ is not of tempered growth and
so the results of [Br] and [BH1] quoted above can not be applied to deduce the
existence of a boundary value.

\example{Example 1.1}
Consider the operator with smooth coefficients 
$$
L=\frac{\partial}{\partial y}-i\frac{2\exp(-y^{-2})}{y^{3}}\frac{\partial}{\partial x}
$$ 
in $\Omega=(-1,1)\times(-1,1)\subset\erre^2$, set
$Z(x,y)=x+i\exp(-y^{-2})$ and define for $y>0$ the function
$f(x,y)=Z^{-1/4}$, where we have used the fact that $\Im Z>0$ for $y>0$
to define the fractional power (we take the branch of $z\mapsto
z^{1/4}$ that is real for $z$ real and positive). For $y>0$ we have
$Lf=0$ and $$ \int_{-1}^1 |f(x,y)|\,dx\le\int_{-1}^1
\frac{1}{\sqrt{|x|}}\le C $$ so, by Theorem 1.1,
$\lim_{y\to0}f(x,y)=bf$ exists and it is easy to check that in fact
$bf(x)=|x|^{-1/2}$.  On the other hand $|f(0,y)|=\exp({y^{-2}/4})$ so
$f$ does not have tempered growth as $y\to0$.

\endexample

\example{Example 1.2}
Let
$$
L=\frac{\partial}{\partial
y}-i\frac{\exp(-y^{-1})}{y^{2}}\frac{\partial}{\partial x} $$ 
in
$\Omega=(-1,1)\times(-1,1)\subset\erre^2$,  set $Z(x,y)=x+i\exp(-y^{-1})$ and
define for $y>0$ the function $f(x,y)=Z^{-1}$. For $y>0$,  $Lf=0$
and $$ \int_{-1}^1 |f(x,y)|\,dx =O(y^{-1}) $$
so, by Theorem 1.1, $\lim_{y\to 0}f(x,y)=bf$ exists (in fact,
$bf(x)=\text{pv}(1/x)-i\pi \delta(x)$). However, $|f(0,y)|=\exp (\frac{1}{y})$.
\endexample

For a holomorphic function $h$ on the rectangle $Q=(-a,a)\times (0,b)$, the
function $h(x+iy)$ is of tempered growth as $y\mapsto 0^+$ if and only if the
integrals 
$$
\int_K|h(x+iy)|dx
$$ 

are of tempered growth. For solutions of a general complex vector
field, the preceding examples indicate that such equivalence is no
longer valid. In the proof of Theorem 1.1, we will use Lemma 1.2 below.
Consider a vector field with smooth coefficients

$$
L=\frac{\partial}{\partial t}+\sum_{j=1}^{n}a_j(x,t)\frac{\partial}{\partial x_j},
$$
defined in a cylinder $D(R,T)=B_R(0)\times(-T,T)\subset\erre^n_x\times\erre_t$, where
$B_R(0)$ denotes the ball $\{x\in\erre^n:\,\,|x|<R\}$. Let $f(x,t)$ and
$g(x,t)$ be two measurable functions in $L^1_{\text{loc}}(D(R,T))$ related by
$$ 
Lf=g\quad\text{in }D(R,T)\tag 1.1
$$
in the sense of distributions.
\proclaim{Lemma 1.2}Let $L$, $f$, $g$ be as above. Then there exist a
continuous function $F(t):(-T,T)\longrightarrow {\D}'(B_R(0))$ and a set
$E\subset(-T,T)$ of Lebesgue measure $|E|=0$ such that 
$$
\<F(t),\psi\>=\int f(x,t)\psi(x)\,dx,\quad t\notin E,\quad
\psi\in\ccinf(B_R(0)). $$
\endproclaim

\demo{Proof of Lemma {\sl 1.2}} 
After shrinking slightly $D(R,T)$ we may assume that $f,g\in L^1(D(R,T))$
and, in view of Fubini's theorem, after modifying $f$ and $g$ on a set of
measure zero  we may also assume that $\int|f(x,t)|\,dx<\infty$ and
$\int|g(x,t)|\,dx<\infty$ for all $|t|\le T$. Fix $\psi\in\ccinf(B_R(0))$. For
any $\phi(t)\in\ccinf(-T,T)$, (1.1) means that  $$
\align
\int_{-T}^T
\bigg(\int f(x,t)\psi(x)\,dx& \bigg)\phi'(t)\,dt=\\
&-\int_{-T}^T
\bigg(\int g\psi-\sum_{j=1}^n\frac{\partial(a_j\psi)}{\partial x_j}
f(x,t)\,dx  \bigg)\phi(t)\,dt.
\endalign
$$
 The expression between parentheses in the right hand side  integral is well defined and defines for each $t$  a distribution $V(t)\in{\D}'(B_R(0))$ of order one such that
$$
\frac{d}{dt}\int f(x,t)\psi(x)\,dx=\<V(t),\psi\>
$$
in the sense of distributions in $(-T,T)$. The function $t\mapsto\<V(t),\psi\>$ is integrable and setting
$$
\<W(t),\psi\>=\int_{0}^t\<V(s),\psi\>\,ds,
$$
it  follows that 
$$
\frac{d}{dt}\bigg(\int f(x,t)\psi(x)\,dx-\<W(t),\psi\>\bigg)=0
$$
in the sense of distributions. Thus, there is a set $E(\psi)$ of measure $|E(\psi)|=0$ such that
$$
\int f(t,x)\psi(x)\,dx-\<W(t),\psi\>=c(\psi),\quad t\notin E(\psi).\tag 1.2
$$
We will now show that $\psi\mapsto c(\psi)$ is a distribution of order one. It is easy to see that $\psi\mapsto c(\psi)$ is linear and if $\{\psi_j\}$ is a sequence converging to zero in $
C^1_{c}(B_R(0))$  then $\<W(t),\psi_j\>\to0$ as $j\to\infty$. 
Using (1.2) for some $t\notin\bigcup_j E(\psi_j)$ we see that $c(\psi_j)\to 0$
when $j\to\infty$. We may now define $F(t)$ by
$$
\<F(t),\psi\>\doteq\<W(t),\psi\>+ c(\psi)= \int_{0}^t\<V(s),\psi\>\,ds+ c(\psi)
$$
and it is clear that the right hand side defines a distribution of order one
in $B_R(0)$ that depends continuously on $t$. Now (1.2) may be restated as 
$$
\int f(x,t)\psi(x)\,dx=\<F(t),\psi\>,\quad t\notin E(\psi).\tag 1.3
$$
We now fix a countable  collection of test functions $\{\psi_j\}$ that is dense in $\ccinf(B_R(0))$ and conclude that (1.3) holds
pointwise for every $\psi\in\ccinf(B_R(0))$ and every $t\notin E=\bigcup
E(\psi_j)$. This proves the lemma. \enddemo

The fact that $F(t)$ is continuous allows us to define the trace of  $f(x,t)$ at $t=t_0$ as $F(t_0)$. This trace 
will in general be just a distribution of order one (a sum of  derivatives
of order $\le1$ of locally finite measures) not representable by a locally
integrable function but, for almost all values of $t$, $F(t)$ is given by the
locally integrable function  $x\mapsto f(x,t)$. 

\example{Example 1.3}
Consider the Mizohata operator $L=\partial_t-it\partial_x$ in
$\Omega=(-1,1)\times(-1,1)\subset\erre^2$ and set $Z=x+it^2/2$,
$f(x,t)=Z^{-1}$. It is easy to check that $f\in L^1(\Omega)$ and that $Lf=0$ in the sense of distributions. For $t\not=0$, $F(t)\in L^\infty(-1,1)\subset L^1(-1,1)$ but for $t=0$ we have $F(0)=\text{pv}(1/x)-i\pi\delta(x)\notin L^1_{\text{loc}}(-1,1)$.
\endexample

The discussion above shows that $f(x,t)$ and $F(t)$ may be identified as distributions in $D(R,T)$. In the sequel, 
we will write just $\int f(x,t)\psi(x)\,dx$ for any value of $t$, even when
the real meaning is $\<F(t),\psi\>$. 

In the next lemma we will need an observation concerning  regularizations of $f$. Let $f(x,t)\in L^1(D(R,T))$ and consider a bump function  $\psi \in C_0^{\infty}(B)$, where $B$ denotes the ball of radius $1$ centered at the origin in ${\Bbb R}^{n+1}$,
of the form $\psi(x,t)=\alpha(x)\beta(t)$.
Assume $\int \alpha(x)\,dx=\int \beta(t)\,dt=1$, and for $\delta >0$, set
$\psi_{\delta}(x,t)=\delta^{-n-1}\psi(x/\delta,t/\delta)=\delta^{-n}\alpha(x/\delta)
\delta^{-1}\beta(t/\delta)=\alpha_\delta(x)\beta_\delta(t)$. 
Extending $f$ as zero outside $D(R,T)$ the convolution $f*\psi_{\delta}(x,t)$ converges to 
$f$ in $L^1(D(R,T))$. Furthermore, for any $\Phi(x,t)\in\ccinf(D(R,T))$ we have
$$
\int f*\psi_{\delta}(x,t)\Phi(x,t)\,dx=
\<F{\buildrel\scriptstyle(t)\over*}\beta_\delta(t),
\Phi(\cdot,t){\buildrel\scriptstyle(x)\over*}\alpha_\delta\>,
$$
where the symbols ${\buildrel\scriptstyle(x)\over*}$ and ${\buildrel\scriptstyle(t)\over*}$ indicate 
convolution in the variables $x$ and $t$ respectively. Since
$\Phi{\buildrel\scriptstyle(x)\over*}\alpha_\delta$ converges in
$\ccinf(B_R(0))$ to $x\mapsto\Phi(t,x)$ uniformly in $t$ and
$F{\buildrel\scriptstyle(t)\over*}\beta_\delta\to F$ uniformly in the
appropriate norm, we may conclude that $$ \lim_{\delta\to0}\int
f*\psi_{\delta}(x,t)\Phi(x,t)\,dx=\<F(t),\Phi(\cdot,t)\> =\int
f(x,t)\Phi(x,t)\,dx $$

\demo{Proof of Theorem {\sl 1.1}} We will use Lemma 1.2 to modify the proof
of Lemma 1.2 in [BH1]. Let $\phi \in C_0^{\infty}(X)$, and $T>0$ such that $$ 
\text{ supp }\phi\times [0,T]\subseteq X\cup U_+ 
$$ 
Shrinking $T$ we may assume without loss of generality that $f$ and $Lf$ are integrable on 
$\text{ supp }\phi\times [\epsilon,T]$ for all $\epsilon>0$. For
$\epsilon \geq 0$ sufficiently small, set
$$
L^{\epsilon}=\frac{\partial}{\partial
t}+\sum_{j=1}^{n}a_j(x,t+\epsilon)\frac{\partial}{\partial x_j}
$$
Let $k\in {\Bbb N}$. We will choose
$\phi_0^{\epsilon},\dots,\phi_k^{\epsilon} \in
C^{\infty}(\overline{U_+})$ such that if 
$$
\Phi^{k,\epsilon}(x,t)=\sum_{j=0}^k\phi_j^{\epsilon}(x,t)\frac{t^j}{j!},
$$ 
then 
$$ 
(1)\quad \Phi^{k,\epsilon}(x,0)=\phi(x),
\quad \text{ and }
\quad (2)\quad |(L^{\epsilon})^{\ast}\Phi^{k,\epsilon}(x,t)|\le Ct^k
$$ 
where $C>0$ depends only  on the size of the derivatives of $\phi$ up to
order $k+1$. In particular, $C$ will be independent of $\epsilon$.  Define
$\phi_0^{\epsilon}(x,t)=\phi(x)$. For $j\geq 1$, write 
$$ 
L^{\epsilon}=\frac{\partial}{\partial
t}+Q^{\epsilon}(x,t,\frac{\partial}{\partial x}) ,
$$ 
and define 
$$
\phi_j^{\epsilon}(x,t)=-\frac{\partial}{\partial t}
\phi_{j-1}^{\epsilon}(x,t)+(Q^{\epsilon})^{\ast}\phi_{j-1}^{\epsilon}
$$
One easily checks that (1) and (2) above hold with these choices
of the $\phi_j^{\epsilon}$. We will next use the integration by
parts formula of the form 
$$
\int u(x,T)w(x,T)\,dx-\int
u(x,0)w(x,0)\,dx=\int_0^T\int_{{\Bbb R}^n}(wPu-uP^{\ast}w)\,dxdt
$$
which is valid for $P$ a vector field, $u$ and $w$ in $C^1({\Bbb
R}^n\times [0,T])$ and the $x-$support of $w$ contained in a
compact set in ${\Bbb R}^n$. Note that the $x$-support of $
\Phi^{k,\epsilon}(x,t)$ is contained in the support of $\phi(x)$.
Let $\psi \in C_0^{\infty}(B_1(0))$, $\psi(x,t)=\alpha(x)\beta(t)$ with  $\int\alpha \,dx=\int\beta dt=1$, as above, and for $\delta >0$, let
$\psi_{\delta}(x,t)=\frac{1}{\delta^{n+1}}\psi(\frac{x}{\delta},\frac{t}{\delta})$. For $\epsilon
>0$, set $f_{\epsilon}(x,t)=f(x,t+\epsilon)$. Observe that if
$\delta < \epsilon$, then the convolution
$f_{\epsilon}*\psi_{\delta}(x,t)$ is $C^{\infty}$ in the region
$t>0$. In the integration by parts formula above set
$u(x,t)=f_{\epsilon}*\psi_{\delta}(x,t)$, $w(x,t)=
\Phi^{k,\epsilon}(x,t)$ and $P=L^{\epsilon}$. We get: 
$$
\align
  \int_Xf_{\epsilon}*\psi_{\delta}(x,0) \phi(x)\,dx &=\int_X f_{\epsilon}*\psi_{\delta}(x,T) \Phi^{k,\epsilon}(x,T)\,dx \\
                 &\qquad-\int_0^T\int_XL^{\epsilon}\left
(f_{\epsilon}*\psi_{\delta}\right )\Phi^{k,\epsilon}\,dxdt \tag 1.4\\            
                   &\qquad\qquad+\int_0^T\int_X
f_{\epsilon}*\psi_{\delta}(L^{\epsilon})^{\ast}\Phi^{k,\epsilon}\,dxdt  
\endalign 
$$ 
Fix $\epsilon >0$. Let $\delta \rightarrow 0^+$. Note that
$f_{\epsilon}*\psi_{\delta}(x,t)$ converges in $L^1$ to $f_{\epsilon}(x,t)$ on
a relatively compact neighborhood $W$ of $\text{ supp } \phi \times [0,T]$.
Hence   
$$ 
L^{\epsilon}\left (f_{\epsilon}*\psi_{\delta}\right )\rightarrow
L^{\epsilon}f_{\epsilon} $$  in ${\D}'(W)$ as $\delta \rightarrow 0^+$.
Moreover, $L^{\epsilon}f_{\epsilon}(x,t)=Lf(x,t+\epsilon)\in L^1$. 
Hence by Friedrichs' Lemma,
$$
L^{\epsilon}\left (f_{\epsilon}*\psi_{\delta}\right )\rightarrow
L^{\epsilon}f_{\epsilon}
$$ 
in $L^1(W)$ as $\delta \rightarrow 0^+$. Finally, the limit as $\delta\to0$
for the first integral on the right hand side of (1.4) was already discussed. 
We thus get  
$$ 
\align
  \int_Xf(x,\epsilon) \phi(x)\,dx &=\int_X f(x,T+\epsilon) \Phi^{k,\epsilon}(x,T)\,dx \\                             
&\qquad-\int_0^T\int_XL^{\epsilon}f_{\epsilon}(x,t)\Phi^{k,\epsilon}(x,t)\,dxdt \tag1.5 \\&\qquad\qquad+\int_0^T\int_X
f_{\epsilon}(x,t)(L^{\epsilon})^{\ast}\Phi^{k,\epsilon}(x,t)\,dxdt \endalign 
$$

  In the third integral on the right, we may integrate first with respect to
$x$, thus obtaining a function of $t\ge0$, $\Gamma_{\epsilon}(t)$
which is  bounded by  
$$ 
|\Gamma_{\epsilon}(t)|
\le \int_X \left|f_{\epsilon}(x,t)(L^{\epsilon})^{\ast}\Phi^{k,\epsilon}(x,t)\right|\, dx 
\leq Ct^k(t+\epsilon)^{-N}\leq Ct^{k-N} ,
$$
 where $C$ depends only on the derivatives of $\phi$ up to order $k+1$ and on the size of its support $K=supp(\phi)$. Furthermore, for $t>0$ we have
$$
\lim_{\epsilon\to0}\Gamma_\epsilon(t)=\lim_{\epsilon\to0}
\int_X f_{\epsilon}(x,t)(L^{\epsilon})^{\ast}\Phi^{k,\epsilon}(x,t)\, dx
=\int_X f(x,t)L^{\ast}\Phi^{k,0}(x,t)\, dx.
$$
 Choose $k=N+1$. By the dominated convergence theorem, as $\epsilon \to 0$,
 this third integral converges to
$\lim_{\epsilon\to0}\int_0^T \Gamma_\epsilon(t)\,dt=\int_0^T\int_X fL^{\ast}\Phi^{k,0}\,dxdt$. 
In the second integral on the right,
note that since $Lf\in L^1(K\times (0,T))$, as $\epsilon\to 0$, the
translates $L^{\epsilon}f_{\epsilon}=(Lf)_{\epsilon}\to Lf$ in $L^1(K\times (0,T))$ while $\Phi^{k,\epsilon} \to\Phi^{k,0}$ uniformly. We thus get  
$$
\langle bf,\phi \rangle = 
\int_X f(x,T)\Phi^{k}(x,T)ds-\int_0^T\int_XLf\Phi^k\,dx dt+
\int_0^T\int_X fL^{\ast}\Phi^k\,dxdt ,
$$ 
where $\Phi^k\doteq\Phi^{k,0}$. From formula
(1.5), we also see that there is $C>0$ independent of $\epsilon$ such that 
$$
|\langle f(.,\epsilon),\phi \rangle |\leq C\sum_{|\alpha |\leq k+1}\Vert
\partial^{\alpha}\phi \Vert_{L^{\infty}}\tag 1.6
$$ 
\enddemo

\proclaim {Corollary 1.3}

Suppose $L=\frac{\partial}{\partial t}+a(x,t)\frac{\partial}{\partial x}$ is
a smooth locally integrable vector field in a neighborhood $U$ of the origin
in the plane. Let $U_+=U\cap {\Bbb R}^2_+$, and suppose $f\in C^0(U_+)$
satisfies $Lf=0$ in $U_+$ and for some integer $N$,
$$\int_K|f(x,t)|dx=O(t^{-N})$$
If the trace $bf=f(x,0)$ is a measure, then it is absolutely continuous with
respect to Lebesgue measure.

\endproclaim

The conclusion in this corollary was proved in [BH1] under the stronger assumptions that
$|f(x,t)|=O(t^{-N})$ and $f\in C^1(U_+)$. However, inspection of the proof
shows that thanks to the formula we have  for the trace $bf$, the proof in
[BH1] goes through with the weaker growth assumption on $f$, at least when
$f\in C^1(U_+)$. To prove it when $f$ is just continuous, we need to make some
modifications in the proof of Theorem 3.1 in [BH1]. Using the notations in
[BH1], we will next indicate the changes needed briefly here and refer the
reader to [BH1] for the details. Recall from [BH1] that for $\zeta$ and $z \in
{\Bbb C}^n$,
$$E(z,\zeta,x,t)=i\zeta\cdot(z-Z(x,t))-\kappa\langle \zeta \rangle
(z-Z(x,t))^2$$
Modifying the notation in [BH1], for $h$ a continuous
function, define
$$g_h(z,\zeta,x,t)=\phi(x)h(x,t)e^{E(z,\zeta,x,t)}$$
where $(z,\zeta)$ are parameters. If $h$ is $C^1$, then for $0<s<t_1$, we have
the analogue of (3.6) in [BH1]:
$$\int_Bg_h(z,\zeta,x,s)d_xZ(x,s)=\int_Bg_h(z,\zeta,x,t_1)d_xZ(x,t_1)+\int_s^{t_1}\int_Bd(g_hdZ)\tag 1.7$$
where $d(g_hdZ)=\left (h(L\phi)+(Lh)\phi\right )e^Edt\wedge dZ$ Suppose now $f$
is continuous and $Lh=0$ for $t>0$. Then if $h_j$ is a sequence of smooth
functions converging uniformly to $f$ in a neighborhood of the closure of
$B\times [s,t_1]$, then since $Lf=0$, by Friedrich's Lemma $Lh_j \mapsto 0$ in
$L^2$ and hence (1.7) will also be valid for $f=h$ leading to
$$\int_Bg(z,\zeta,x,s)d_xZ(x,s)=\int_Bg(z,\zeta,x,t_1)d_xZ(x,t_1)+\int_s^{t_1}\int_Bd(gdZ)\tag 1.8$$
where 
$$g(z,\zeta,x,t)=\phi(x)f(x,t)e^{E(z,\zeta,x,t)}$$
and $$d(gdZ)=fL\phi e^Edt\wedge dZ$$
Suppose now the integrals $\int_K|f(x,t)|dx$ have a tempered growth as in
Theorem 1.1. Then as $s\mapsto 0$, the integral on the left in (1.8) converges
to $\int_Bg(z,\zeta,x,0)d_xZ(x,0)$. We claim that for the directions $\zeta$ of
interest in Theorem 3.1, as $s\mapsto 0$, the second integral on the right in
(1.8) converges to 

$$\int_0^{t_1}\int_Bd(gdZ)$$
Indeed, the latter follows from Lemma 1.2 which tells us that the
distributions $f(.,t)$ are uniformly bounded which in our
situation implies an exponential decay in the $x$ integral. One can then use
the Dominated Convergence Theorem to prove the assertion.

\heading{2. On pointwise convergence of solutions to their traces}
\endheading

Suppose  $L$ is a never vanishing smooth vector field on a domain D in
the plane and $f$ is a smooth solution of $Lf=0$ in $D$ with tempered growth as
one approaches a noncharacteristic  boundary piece $\gamma$. Assume that on
$\gamma$ the function $f$ has a weak boundary value $bf$ which is locally
integrable. Unlike the case of the Cauchy Riemann operator, simple examples
show that even when $L$ is real analytic, $f$ may not converge nontangentially
to $bf$. Indeed, consider the Mizohata vector field
 $$
L_1=\frac{\partial}{\partial t}-2it\frac{\partial}{\partial x}
$$
 Let $F(z)$ be a holomorphic function in the semi-disc
$\{z=x+iy:|z|<1, y>0\}$ with a weak $L^1$ boundary value $bF$ on the $x$-axis. If $F$ is chosen so that it is bounded and  on a set of full measure 
in $(-1, 1)$ it has no limit in parabolic regions, then  the function
$F(x+it^2)$ is a solution of $L$ which does not converge nontangentially to
its weak limit $bF$ as $t$ tends to $0$. The existence of such $F$ follows
from the following more general theorem in [Z].

\proclaim {Theorem 7.44} ([Z]) Let $C_0$ be any simple closed curve passing
through $z=1$, situated, except for that point, totally inside the circle
$|z|=1$, and tangent to the circle at that point. Let $C_{\theta}$ be the
curve $C_0$ rotated around $z=0$ by an angle $\theta$. There is a Blaschke
product $B(z)$ which, for almost all $\theta_0$, does not tend to any limit as
$z\mapsto \exp (i\theta_0)$ inside $C_{\theta_0}$.
\endproclaim

The preceding theorem shows that even for the $C^{\infty}$ and analytic
hypoelliptic vector field 
$$L_2=\frac{\partial}{\partial t}-3it^2\frac{\partial}{\partial x}$$
we can get bounded solutions $f(x,t)=F(x+it^3)$ with $F$ holomorphic in a
semi-disc in the upper half plane, $bf\in L^1$ but $f(x,t)$ does not converge
nontangentially on a subset of full measure in $(-1,1)$. For both examples
$L_1$ and $L_2$, the solutions $f(x,t)$ converge to their boundary values a.e.
in certain cuspidate regions containing the vertical segments $\{(x,t):t>0\}$.
However, if we consider the vector field $L_3=\frac{\partial}{\partial t}$,
and take $f(x,t)=bf(x)=$ the characteristic function of a Cantor set $C$  of positive
measure in an interval $(a,b)$, the only sets of approach for which $f(x,t)\to bf(x)$ for a.e. $x\in C$  are the vertical segments. Therefore, for a general locally
integrable vector field, we can not get sets of approach for convergence
larger than curves. In this section we prove a.e. convergence along vertical
segments to $L^1$ boundary values for locally integrable vector fields of the
form
$$L=\frac{\partial}{\partial t}+a(x,t)\frac{\partial}{\partial x}$$
In section 5 we will prove that when $L$ is a locally solvable vector field,
at appropriate points, the sets of approach for a.e. convergence are open
sets.\newline To state our result in an invariant form, let $\Omega$ be a
smooth domain in the plane, $L=X+iY$ a locally integrable vector field near
each point of a piece $\Sigma$ of the boundary and $f$ a continuous solution
of $Lf=0$ in $\Omega$. Assume that for some defining function $\rho$ of
$\Omega$, there exists an integer $N$ such that the line integrals
$$ \int_{\rho=t}|f|d\sigma_t=O(t^{-N})$$
Suppose $\Sigma$ is noncharacteristic for $L$. Then by Theorem 1.1, $f$ has a
trace $bf$ on $\Sigma$. Assume that $bf\in L^1(\Sigma)$. After contracting
$\Sigma$ around one of its points, we can get a smooth first integral $Z$ for
$L$ with the property that the restriction of $\Re Z$ on $\Sigma$ has a nonzero
differential. For $p\in \Sigma$, the set
$$S(p)=\{w:\Re Z(w)=\Re Z(p)\}$$
is a curve near $p$ which is transversal to $\Sigma$. Let $S^+(p)$ denote the
part of this curve lying in $\Omega$. We will prove:

\proclaim {Theorem 2.1} For almost all $p\in \Sigma$,
$$\lim_{S^+(p)\owns q\mapsto p}f(q)=bf(p)$$
\endproclaim

To prove Theorem 2.1, we first flatten $\Sigma$ in new coordinates. By
hypotheses, $x=\Re Z$ and $t=\rho$ form a change of coordinates near a fixed
point $p\in \Sigma$ which we may assume is mapped to the origin. In these
coordinates, $\Sigma$ and $L$ are given by 
$$t=0,\qquad L=\lambda(x,t)\left (\frac{\partial}{\partial
t}+a(x,t)\frac{\partial}{\partial x}\right )$$
respectively, for some nonvanishing smooth factor $\lambda$, and the first
integral takes the form
$$Z(x,t)=x+i\varphi(x,t)$$
for some real-valued $\varphi$. Therefore, Theorem 2.1 follows from

\proclaim{Theorem 2.1'} Let
$$
L=\frac{\partial}{\partial t}+a(x,t)\frac{\partial}{\partial x}
$$
be a smooth locally integrable complex vector field in the subset
$U=(-r,r)\times (-T,T)$ of the plane. Assume $f$ is continuous on $U^+=(-r,r)\times (0,T)$ and 
$Lf=0$ in $U^+$. Suppose for any  compact set $K\subset\subset (-r,r)$ there
exists $N=N(K)\in {\Bbb N}$, and $C=C(K)>0$ such that 
$$
  \int_K |f(x,t)|\,dx \le \frac{C}{t^N} \quad, \text{ as } t\to 0^+
 $$
and the boundary value $bf\in L^1(-r,r)$. Then 
$$
\lim_{t\to 0}f(x,t)=bf(x)
$$
a.e. in $(-r, r).$
\endproclaim

We begin with some general lemmas which are valid for a
general, not necessarily locally integrable vector field.

\proclaim {Lemma 2.2} Let 
$$ 
L={\partial\over\partial
t}+\sum_{j=1}^n b_j(x,t){\partial\over\partial x_j} 
$$ 
be smooth
on a neighborhood $U=B(0,a)\times(-T,T)$ of the origin in
$\erre^{n+1}$ with $B(0,a)=\{x\in\erre^n:\,\,\,|x|<a\}$. We will
assume that the coefficients $b_j(x,t)$, $j=1,\dots,n$ vanish on $F\times[0,T)$,  where $F\subset
B(0,a)$ is a closed set.   Assume that  $f$  is continuous  on
$U^+=B(0,a)\times(0,T)$, satisfies $Lf=0$ in $U^+$ and 
for any  compact set $K\subset\subset B(0,a)$, there exists $N=N(K)\in {\Bbb
N}$, and $C=C(K)>0$ such that 
$$
  \int_K |f(x,t)|\,dx \le \frac{C}{t^N} \quad, \text{ as } t\to 0^+
 $$
 and $bf\in L^{1}(B(0,a))$. Then pointwise, 
$$
\lim_{t\to 0}f(x,t)=bf(x)\quad
\text{a.e. on}\quad F
$$ 
\endproclaim

\proclaim {Remark} The existence of a distribution boundary value $bf$ was
proved in Theorem 1.1.
\endproclaim
In the proof of this lemma, we will use another lemma which we will first
prove. \proclaim{Lemma 2.3} Let
$$ 
L={\partial\over\partial t}+\sum_{j=1}^n a_j(x,t){\partial\over\partial x_j}
$$
be a smooth complex vector field on an open set $U=B(0,r)\times (0,T)$ in
${\Bbb R}^{n+1}$ . Assume $f$ is continuous in $U$ and satisfies $Lf=0$ in
$U$. Suppose $a_j(0,t)=0$ for all $j$ and  for all $t\in (0,T)$. Then $f(0,t)$
is constant.
\endproclaim

\demo{Proof}
 Let $\phi(x)\in C_0^{\infty}(B(0,1))$ such that the sequence
$\phi_{\epsilon}(x)=\frac{1}{{\epsilon}^n}\phi(\frac{x}{\epsilon})$ forms an
approximate identity family. Using $Lf=0$ and integration by parts, for any
$0<a<b<T$, we have:
$$
\int f(x,b)\phi_{\epsilon}(x)\,dx-\int  
f(x,a)\phi_{\epsilon}(x)\,dx=-\int_a^b\int_{{\Bbb
R}^n}f(x,t)L^t\phi_{\epsilon}(x)\,dxdt \tag 2.1
$$
Observe that the left hand side converges to $f(0,b)-f(0,a)$ as $\epsilon \to
0$. It therefore suffices to show that the right hand side converges to $0$.
We write
$$
\aligned
\int_a^b\int_{{\Bbb R}^n}f(x,t)L^t\phi_{\epsilon}(x)\,dxdt
&=\int_a^b\int_{{\Bbb
R}^n}(f(x,t)-f(0,t))L^t\phi_{\epsilon}(x)\,dxdt\\
&+\int_a^b\int_{{\Bbb
R}^n}f(0,t)L^t\phi_{\epsilon}(x)\,dxdt
\endaligned
$$
Note that
$$
\int_a^b\int_{{\Bbb R}^n}f(0,t)L^t\phi_{\epsilon}(x)\,dxdt=-\int_a^bf(0,t)\left(\sum_{j=1}^n\int_{{\Bbb R}^n}\frac{\partial}{\partial x_j}(a_j\phi_{\epsilon}(x))\,dx\right)dt=0
$$
since $\phi_{\epsilon}(x)$ has compact support. We now estimate the other term
$$
\align
\left |\int_a^b\int_{{\Bbb R}^n}(f(x,t)-f(0,t))
L^t\phi_{\epsilon}(x)\,dxdt\right |&\leq \\
\Bigg|\int_a^b\int_{{\Bbb R}^n}(f(x,t)-f(0,t))\,\text{div\,}(a)&\phi_{\epsilon}(x)\,dxdt\Bigg|\\
+\Bigg |\sum_j\int_a^b&\int_{|x|\leq \epsilon}(f(x,t)-f(0,t))a_j\frac{\partial
{\phi_{\epsilon}}}{\partial x_j}\,dxdt\Bigg|\\
&\leq
  C\omega(\epsilon)+C\omega(\epsilon)\frac{1}{{\epsilon}^{n+1}}\int_{|x|\leq
  \epsilon} |x|\,dx\\
&\leq C_1\omega(\epsilon) 
\endalign
$$
where in the second inequality we have used the vanishing of the $a_j(0,t)$
and used the notation  $\omega(\epsilon)=\text{ sup }|f(x,t)-f(0,t)|$ on
$B(0,\epsilon)\times [a, b]$. Since $f$ is continuous, it follows that
$\omega(\epsilon)\to 0$ and hence $f(0,t)$ is constant.
\enddemo
\demo{Proof of Lemma \sl 2.2} By Lemma 2.3, for any $x\in F$,
$f(x,t)=f(x,T)$. We therefore have to show that $bf(x)=f(x,T)$ a.e. in $F$. We
recall from [BH1] (see the proof of Lemma 3.3) that for any  $\phi\in
C_c^{\infty}(B(0,a))$, and any $k\in \Bbb N$, we can choose smooth functions
$\phi_0,\dots,\phi_k$ with the properties that if  
$$
\Phi^j(x,t)=\sum_{l=0}^j\phi_l(x,t)\frac{t^l}{l!}\quad \text{for} \quad j\leq k 
$$  
then
$$
(1) \quad \Phi^j(x,0)=\phi(x),\quad \text{and } (2)\quad 
L^t\Phi^j(x,t)=\frac{\phi_{j+1}}{j!}t^j
$$
Moreover, since the coefficients $b_j(x,t)$ vanish on $F\times [0,T]$, each
$\phi_j$ has the form 
$$
\phi_j(x,t)=\sum_{|\alpha|\leq j}c_{\alpha}(x,t)D_x^{\alpha}\phi(x)
$$
where the $c_{\alpha}$ are smooth and satisfy the estimate
$$ 
|c_{\alpha}(x,t)|\leq Cd(x,F)^{|\alpha|}\tag 2.2
$$
where $d(x,F)$ denotes the distance from $x$ to $F$ . The constant  $C$ in
(2.2) is independent of the $\phi_j$ since the $c_{\alpha}$ are obtained from
the coefficients $b_j(x,t)$ of $L$ by means of algebraic operations and
differentiations. The proof of Theorem 1.1 also shows us that  
$$
\langle bf,\phi
\rangle=\int_{B(0,a)}f(x,s)\Phi^k(x,s)\,dx+\int_0^s\int_{B(0,a)}f(x,t)L^t\Phi^k(x,t)\,dxdt\tag 2.3
$$
if we choose $k=N+1$. Let $K\subseteq F$ be a compact set. Choose smooth
functions $0\leq \phi_{\epsilon}(x)\leq 1$ in $C_c^{\infty}(B(0,a))$
satisfying: $(1)\quad \phi_{\epsilon}(x)=1$ for $x\in K$; $(2)
\quad \phi_{\epsilon}(x)=0$ if $d(x,K)>\epsilon$; and $(3)
\quad |D_x^{\alpha}\phi_{\epsilon}(x)|\leq A_{\alpha}\epsilon^{-|\alpha|}$.
Thus $\phi_{\epsilon}(x)$ converges pointwise to the characteristic function
of $K$ and for $|\alpha|>0$, $D^{\alpha}\phi_{\epsilon}(x)\to 0$. Let $\psi\in
C_c^{\infty}(B(0,a))$ and apply (2.3) to $\phi=\phi_{\epsilon}\psi$ to get $$
\align
\langle
 bf,\phi_{\epsilon}\psi\rangle&=\int_{B(0,a)}f(x,s)\Phi^{k,\epsilon}(x,s)
\,dx\\
&+\int_0^s\int_{B(0,a)}f(x,t)L^t\Phi^{k,\epsilon}(x,t)\,dxdt\tag 2.4
\endalign
$$ 
Since the sequence $\phi_{\epsilon}\psi$ is uniformly bounded , converges
pointwise to $\psi(x)\chi_{K}(x)$, and $bf$ is integrable, by the dominated
convergence theorem,  
$$
\langle bf,\phi_{\epsilon}\psi\rangle\to \int_Kbf(x)\psi(x)\,dx\quad \text{
as } \epsilon \to 0
$$ 
We consider next the first integral on the right in
(2.4): 
$$
\int_{B(0,a)}f(x,s)\Phi^{k,\epsilon}(x,s)\,dx=\sum_{j=0}^k\int_{B(0,a)}f(x,s)\phi^{\epsilon}_j(x,s)\frac{s^j}{j!}\,dx\tag 2.5$$ Recall that  
$$
\phi^{\epsilon}_j(x,s)=\sum_{|\alpha|\leq
j}c_{\alpha}(x,s)D_x^{\alpha}(\phi_{\epsilon}(x)\psi(x))
=\sum_{|\alpha|\leq j}\sum_{\beta\leq
\alpha}c_{\alpha,\beta}(x,s)D_x^{\beta}\phi_{\epsilon}(x)D_x^{\alpha-\beta}\psi(x)\tag 2.6$$
In the double sum above, if $\beta<\alpha$, then
$$
\align
|c_{\alpha,\beta}(x,s)D_x^{\beta}\phi_{\epsilon}(x)D_x^{\alpha-\beta}\psi(x)|
&\leq Cd(x,F)^{|\alpha|}|D_x^{\beta}\phi_{\epsilon}(x)|\\
&\leq Cd(x,K)^{|\alpha|}|D_x^{\beta}\phi_{\epsilon}(x)|\\
&\leq C{\epsilon}^{|\alpha|-|\beta|}\tag 2.7
\endalign
$$
Hence such terms go to $0$ as $\epsilon \to 0$. Therefore, we only need to
look at the contribution of 
$$
\sum_{|\alpha|\leq j}c_{\alpha}(x,s)(D_x^{\alpha}\phi_{\epsilon}(x))\psi(x)
$$
In this latter sum, when $|\alpha|\geq 1$, the term
$$
c_{\alpha}(x,s)(D_x^{\alpha}\phi_{\epsilon}(x))\psi(x) \to 0
$$
pointwise and the sequence is uniformly bounded independently of $\epsilon$.
Therefore, by the dominated convergence theorem, 
$$
\int c_{\alpha}(x,s)D_x^{\alpha}\phi_{\epsilon}(x)\psi(x)f(x,s)\,dx \to 0\quad
\text{as } \epsilon \to 0
$$
It follows that when $j\geq 1$,
$$
\align
\lim_{\epsilon \to 0}\int_{B(0,a)}f(x,s)\phi_j^{\epsilon}(x,s)\frac{s^j}{j!}\,dx
&=\lim_{\epsilon \to
0}\int_{B(0,a)}f(x,s)\phi_{\epsilon}(x)\psi(x)c_0(x,s)s^j\,dx\\
&=\left(\int_Kf(x,s)\psi(x)c_0(x,s)\,dx\right)s^j\\
&=\left(\int_Kf(x,T)\psi(x)c_0(x,s)\,dx\right)s^j\tag 2.8
\endalign
$$ 
where we used Lemma 2.3 in the last equation. Since
$\phi_0^{\epsilon}(x,s)=\phi_{\epsilon}(x)\psi(x)$, from (2.5) and (2.8) we
see that  
$$
\align
\lim_{\epsilon\to
0}\int_{B(0,a)}f(x,s)\Phi^{k,\epsilon}(x,s)\,dx&=\lim_{\epsilon\to
0}\int_{B(0,a)}f(x,T)\phi_{\epsilon}(x)\psi(x)\,dx+O(s)
\\
&=\int_Kf(x,T)\psi(x)\,dx+O(s)\tag 2.9
\endalign
$$
Consider next the double integral in (2.4). Since
$$
L^t\Phi^{k,\epsilon}(x,t)=\frac{\phi^{\epsilon}_{k+1}t^k}{k!},\quad k=N+1
$$
and  $\phi^{\epsilon}_{k+1}$ is bounded independently of $\epsilon$, it follows that
$$
\int_0^s\int_{B(0,a)}f(x,t)L^t\Phi^{k,\epsilon}(x,t)\,dxdt=O(s^2)\tag 2.10
$$
Finally from (2.4), (2.9), and (2.10), we get:
$$
\int_Kbf(x)\psi(x)\,dx=\int_Kf(x,T)\psi(x)\,dx+O(s)
$$
Letting $s\to 0$ in the latter, we conclude that $bf(x)=f(x,T)$ a.e. in $K$
and hence in $F$.
\enddemo

\demo{Proof of Theorem \sl 2.1'} We may assume that $L$ has a smooth first
integral $Z$ with the property that $LZ=0$ in $U$ and the differential
$dZ(x,t)\ne 0$ at every point in $U$. We may in fact assume that 
$Z(x,t)=x+i\phi(x,t)$, where $\phi$ is real valued, $\phi(0,0)=0$,
$D_x\phi(0,0)=0$ and $D_x^2\phi(0,0)=0$.  Let
$$
E=\{x\in (-r,r):\exists \quad \epsilon >0\text{ with } \phi(x,t)\equiv
\phi(x,0)\quad \forall t\in [0,\epsilon]\}
$$ 
Then by Lemma 2.2, for almost all points
in $E$, $\lim_{t\to 0}f(x,t)=bf(x)$. Consider therefore a point $x_0\notin E$,
say $x_0=0\notin E$. Then there exists a sequence $t_j$ decreasing to zero
such that $\phi(0,t_j)\neq 0$. After decreasing $r$ and $T$, by the boundary
analogue of the Baouendi-Treves approximation theorem (see Theorem 3.1 in
[T2] ), there is a sequence of entire functions $P_k$ such that
$P_k(Z(x,t)) \to f(x,t)$ in the sense of distributions on $U^+=(-r,r)\times
(0,T)$. If there are two sequences $\{s_j\}$ and $\{y_k\}$ both converging to
$0$ with $\phi (0,s_j)>0$ and $\phi (0,y_k)<0$, then the image $Z(U^+)$ will
contain a ball $B$ centered at $Z(0,0)=0$ on which the  entire functions $P_k$
will converge uniformly to a holomorphic function $H$ and so $f(x,t)=H(Z(x,t))$
will in fact be smooth up to $t=0$. Indeed, this latter assertion follows for a $C^1$ $f$ since 
Theorem 3.1 in [T2]  guarantees uniform convergence on
compact subsets of $U^+=(-r,r)\times (0,T)$ for such solutions. In the general case, 
we can use the representation formula of Theorem 6.4 in [T2] to express $f$
as $Qh$ where $h$ is a $C^1$ solution and $Q$ is a second order elliptic
differential operator which maps solutions to solutions. We can then get a
holomorphic function $G$ on the ball $B$ such that $ h(x,t)=G(Z(x,t))$ and  so
from the form of the operator $Q$, $f$ will also equal $P(Z(x,t))$ for some
holomorphic function $P$ on $B$. We may therefore assume that $\phi (0,t)$
does not change sign on some interval $[0,T]$. Without loss of generality, we
may assume that  
$$
\phi (0,t)\geq 0 \text{ for }t\in [0,T]\tag 2.11
$$ 
We also have a sequence $t_j$ converging to $0$ where now $\phi(0,t_j)>0$ for all $j$.  By
Theorem 3.1 in [BH1], it follows that at the origin, the FBI transform  (with
$Z(x,0)$ as phase) of  $bf(x)$ decays exponentially in a complex conic
neighborhood of the covector $(0;-1)$. By Theorem 2.2 in [BCT], there exists
an interval centered at the origin which we will continue to denote by 
$(-r,r)$, a number $\delta >0$ and a holomorphic function $F$ of tempered
growth defined on the open set 
$$ 
Q=\{Z(x,0)+iZ_x(x,0)v:x\in (-r,r),0<v<\delta\}
$$ 
such that for any $\psi\in C_c^{\infty}(-r,r)$, 
$$ 
\int
bf(x)\psi(x)\,dx=\lim_{v\to 0}\int F(Z(x,0)+iZ_x(x,0)v)\psi(x)\,dx
$$ 
Since $bf$ is a locally integrable function, as is well known, the holomorphic function $F$
converges nontangentially to $bf(x)$ almost everywhere (see for example,
Corollary 1.1 in [BH2]). We may assume that $0$ is a point where this
convergence holds.  Let $M$ be a vector field for which the function
$Z_1(x,t)=Z(x,0)+iZ_x(x,0)t$ is a first integral. Note that the function
$F(Z_1(x,t))$ is a solution of $M$ for $x$ near $0$ and $t>0$. We can
therefore apply the boundary version of the approximation theorem both to $M$
with the  solution  $F(Z_1(x,t))$ and to $L$ with the solution $f(x,t)$ to
deduce the following: for $(x,t)\in (-a,a)\times (0,r)$, $a $ and $r$
sufficiently small, in the distribution sense :    
$$ 
f(x,t)=\lim_{\tau\to\infty}\left(\frac{\tau}{\pi}\right
)^{\frac{1}{2}}\int_We^{-\tau(Z(x,t)-Z(y,0))^2}bf(x)g(y)\,dZ(y,0)  
$$   
and likewise  
$$ 
F(Z_1(x,t))=\lim_{\tau\to\infty}\left
(\frac{\tau}{\pi}\right)^{\frac{1}{2}}
\int_We^{-\tau(Z_1(x,t)-Z(y,0))^2}bf(x)g(y)\,dZ(y,0)  
$$  
where $g$ is a smooth function supported in some neighborhood of $0$, identically equal to
$1$ near $0$. In the above limits, we have taken advantage of the fact that
$Z(x,0)\equiv Z_1(x,0)$. We observe that the second limit is valid since the
$x$ derivative of $Z_1(x,0)$ at $0$ is $1$ (see Theorem 3.1 in [T2]). These
formulas show that there exist entire functions $P_{\tau}(z)$ such that in the
distribution sense,  
$$  
f(x,t)=\lim _{\tau \to \infty}P_{\tau}(Z(x,t))\quad
\text{ and } \quad F(Z_1(x,t))=\lim _{\tau \to \infty}P_{\tau}(Z_1(x,t))  
$$  
We observe that since the vector field $M$ is elliptic near the origin,
$P_{\tau}(Z_1(x,t))$ converges uniformly on compact subsets of $(-a,a)\times
(0,r)$ to $F(Z_1(x,t))$. Recall now that  
$$ 
\lim_{t\to 0}F(Z_1(0,t))=\lim_{t\to 0}F(it)=bf(0)
$$ 
Let $\epsilon >0$. There exists $\delta >0$ such that if $0<t<\delta$ and $\phi (0,t)>0$, then  
$$ 
|F(i\phi (0,t))-bf(0)|<\epsilon
$$ 
Let $0<s<\delta$. We consider two cases on $\phi_t(0,s)$. Assume first that  $\phi_t(0,s)\neq 0$. Then the vector field $L$ is elliptic at $(0,s)$ and hence the sequence $P_{\tau}(Z(x,t))$ converges uniformly to $f(x,t)$ near $(0,s)$. In particular,  
$$
\lim _{\tau \to \infty}P_{\tau}(Z(0,s))=f(0,s)
$$ 
If in addition, $\phi (0,s)>0$, then 
$$
f(0,s)=\lim _{\tau \to \infty}P_{\tau}(Z(0,s))=\lim _{\tau \to
\infty}P_{\tau}(Z_1(0,\phi (0,s)))=F(i\phi (0,s)) 
$$ 
Hence,  
$$
|f(0,s)-bf(0)|<\epsilon 
$$  
If $\phi (0,s)=0$, then since $\phi _t(0,s)\neq 0$, there exists $y$ arbitrarily close to $s$ where $\phi (0,y) >0$ and $\phi_t(0,y)\neq 0$ and so as we already saw, we will still have
$|f(0,s)-bf(0)|<\epsilon$. Suppose now $\phi _t(0,s)=0$. If $s$ is in the
closure of  
$$
\{y:0<y<\delta, \phi_t(0,y)\neq 0\}
$$ 
then by the first case and continuity of $f$, $|f(0,s)-bf(0)|<\epsilon$. If $s$ is not in the closure of this set, then since there is a 
sequence $t_j\to 0$ where $\phi(0,t_j)>0$ and
$\phi (0,0)=0$, we can find $y\in (0,s)$ such that $\phi _t(0,t)=0$ on the
interval $(y,s)$, and $\phi_t(0,y)>0$. By Lemma 2.3, we will then have
$f(0,s)=f(0,y)$. Hence, in this case too we get: 
$$ 
|f(0,s)-bf(0)| =|f(0,y)-bf(0)|<\epsilon 
$$
 
\enddemo

\heading {3. Locally solvable vector fields and Hardy spaces}
\endheading

Consider a locally solvable vector field with smooth coefficients 
 $$
L=\frac{\partial}{\partial y}+a(x,y)\frac{\partial}{\partial x}
$$
defined on a neighborhood of the origin $Q=[-a,a]\times[-b,b]$. Since our point of view is local and locally 
solvable vector fields are known to be locally integrable [T1], we will
assume without loss of generality that there is a smooth real function
$\varphi(x,y)$ defined on a neighborhood of $Q$ such that
$Z(x,y)=x+i\varphi(x,y)$ is a first integral of $L$, i.e., $LZ=0$ or,
equivalently, $a(x,y)=-i\varphi_y(x,y)/(1+i \varphi_x(x,y))$. Furthermore, it
is convenient for technical reasons to assume as well that
$\varphi(0,0)=\varphi_x(0,0)=0$ and 
$$ 
|\varphi_x(x,y)|<\frac{1}{2}\quad\text{on a neighborhood of $Q$.}\tag 3.1  
$$ 
It is well known that the local solvability of $L$ is equivalent to the  fact that $L$ satisfies the Nirenberg-Treves condition $(\P)$ ([NT],[T1]) and this reflects on the
behavior of $\varphi$ in the following way: 
$$ 
\text{\it for every }x\in[-a,a]
\text{\it\ the map } [-b,b]\owns y\mapsto\varphi(x,y) \text{\it\ is
monotone.}
$$ 
We can now state the main result of this paper: 
\proclaim {Theorem 3.1} Suppose $f$ is continuous and is a weak solution of
$Lf=0$ in the rectangle $(-a^{\prime},a^{\prime})\times (0,b^{\prime})$ for
some $a^{\prime}>a$, $b^{\prime}>b$. Assume that there is a positive integer
$N$ such that for each $K$ compact in $(-a^{\prime},a^{\prime})$,
$\int_K|f(x,y)|\,dx=O(y^{-N})$. Suppose the boundary value of $f$ at $y=0$,
$bf\in L^p(-a^{\prime},a^{\prime})$ for some $1\leq p \leq \infty$. Then for
any $a<a^{\prime}$, the norms of the traces $f(.,y)$ in $L^p[-a,a]$ are
uniformly bounded as $y\mapsto 0^+$. \endproclaim The proof of Theorem 3.1
will occupy most of the rest of the paper. We begin by defining 
$$ 
m(x)=\min_{0\le y\le b} \varphi(x,y)\qquad  M(x)=\max_{0\le y\le b} \varphi(x,y),\quad -a\le x\le a. 
$$
Thus, the function $Z(x,y)$ takes the rectangle $Q=[-a,a]\times [0,b]$ onto 
$$
Z(Q)=\{\xi+i\eta:\quad -a\le\xi\le a,\quad m(\xi)\le\eta\le M(\xi)\}.
$$
The interior of $Z(Q)$ is 
$$
\{\xi+i\eta:\quad -a<\xi<a,\quad m(\xi)<\eta< M(\xi)\},
$$
in particular, this interior is not empty if and only if $M(x)>m(x)$ for some $x\in(-a,a)$. 
The case of an empty interior corresponds to the uninteresting and trivial
case in which $\varphi$ is independent of $y$ and $L=\partial_y$ so we will
assume from now on that $\varphi_y$ does not vanish identically for $y>0$
which in particular implies that  $Z(Q^+)$ has nonempty interior. Every
connected component $U$ of the interior of $Z(Q^+)$ is of the form 
$$
U=\{\xi+i\eta:\quad \alpha<\xi<\beta,\quad m(\xi)<\eta< M(\xi)\}, 
$$ 
where $(\alpha,\beta)$ is a connected component of the open set $\{x\in(-a,a):\quad
M(x)>m(x)\}$. Notice that, by the very definition of $U$, it follows that
$M(\alpha)=m(\alpha)$ unless $\alpha=-a$, and $M(\beta)=m(\beta)$ unless
$\beta=a$. We will focus our attention on the case where  $-a<\alpha<\beta<a$,
so $M(\alpha)=m(\alpha)$ and $M(\beta)=m(\beta)$. Notice that, because for
every $x\in(\alpha,\beta)$ the map $y\mapsto\varphi(x,y)$ is monotone and not
constant, it is clear that  either $\varphi_y(x,y)\ge0$ for all
$x\in(\alpha,\beta)$ and $|y|\le b$ or $\varphi_y(x,y)\le0$ for all
$x\in(\alpha,\beta)$ and $|y|\le b$. From now on we will assume that the first
possibility occurs, i.e., that $\varphi_y\ge0$ on
$[\alpha,\beta]\times[-b,b]$. Hence,  
$$ 
M(x)=\varphi(x,b)\quad \text{ and }\quad m(x)=\varphi(x,0), \quad \alpha\le x\le\beta. 
$$ 
Thus, $U$ is a bounded region lying between two smooth graphs and its boundary $\partial
U$ is smooth except at the two end points $(\alpha,M(\alpha))$ and
$(\beta,M(\beta))$. Note that  $U$ has a rectifiable boundary
of length bounded by 
$$
\align 
|\partial U|&\le\int_\alpha^\beta\sqrt{1+\varphi^2_x(x,b)}\,dx + 
\int_\alpha^\beta\sqrt{1+\varphi^2_x(x,0)(x)}\,dx\\ &\le 2(\beta-\alpha)
\sqrt{1+\sup_Q|\nabla\varphi|^2} = K(\beta-\alpha). 
\endalign 
$$

We will first show that our solution $f$ determines a holomorphic function
$F$ on $U$ such that $f(x,y)=F(Z(x,y))$. In order to see this we set, for small
$\epsilon>0$ and big $\tau>0$, 
$$ 
E_{\tau,\epsilon} f(x,y)=
({\tau/\pi})^{1/2}\int_\erre e^{-\tau
 [Z(x,y)-Z(x',\epsilon)]^2} f(x',\epsilon) h(x') Z_x(x',0)\, dx'. \tag 3.2
$$
Here $h(x')\in\ccinf(-a',a')$ is a test function identically equal to 1 on a
neighborhood of $[-a,a]$. Thanks to assumption (3.1), the proof of the
Baouendi-Treves approximation theorem [BT] implies that, for fixed $\epsilon$,
$E_{\tau,\epsilon} f(x,y)\to f(x,y)$ uniformly on the rectangle $R_\epsilon$
given by $|x|\le a$, $\epsilon\le y\le b$, provided $b$ is small enough (to be
more specific, provided $b\,\sup|\nabla\varphi|<<1$) which we could have
 assumed from the start. Formula (3.2) may be written as $E_{\tau,\epsilon}
f(x,y)=F_{\tau,\epsilon}(Z(x,y))$ where $F_{\tau,\epsilon}$ is an entire
function. If we take  a sequence $\tau_k\to\infty$, we conclude that
$F_{\tau_k,\epsilon}$ is uniformly Cauchy on compact subsets of 
the set $Z(R_\epsilon)$. In particular, $F_{\tau_k,\epsilon}$ converges
uniformly to a  function $F_{\epsilon}$ which is holomorphic on 
$$
U_\epsilon=\{(\xi+i\eta):\quad \alpha<\xi<\beta,\quad \varphi(\xi,\epsilon)<
\eta< M(\xi)\} 
$$ 
and continuous on  
$$ 
\{(\xi+i\eta):\quad\alpha<\xi<\beta,\quad \varphi(\xi,\epsilon)\le \eta\le M(\xi)\}. 
$$ 
Thus,
$F_\epsilon(Z(x,y))=f(x,y)$ on $\alpha<x<\beta$, $\epsilon\le y<b$, and
$F_\epsilon$ is an extension of $F_{\epsilon'}$ if $0<\epsilon<\epsilon'$ are
small. As $\epsilon\searrow0$ we obtain a holomorphic function $F$ defined on
$Z(U)$ such that $f(x,y)=F(Z(x,y))$. We will study the boundary limits of
$F$ in $U$. Since $F$ is continuous on the graph $\Gamma$ of $M(x)$,
$\alpha<x<\beta$, it is apparent that the boundary value $bF$ of $F$ in
$\Gamma$ is given by $bF(x+iM(x))=bF(x+i\varphi(x,b))=f(x,b)$ and we need only
worry about the behavior of $F$ when approaching the lower part $\gamma$ of
$\partial U$ given by the  graph $\eta=\varphi(\xi,0)$. We now recall the
definition of a Hardy space (see [Du]) for a domain with rectifiable boundary.
 
\proclaim {Definition 3.1}For $1\leq p<\infty$, a holomorphic function $g$ on a
bounded domain $D$ with rectifiable boundary is said to be in $E^p(D)$ if
there exists a sequence of rectifiable curves $C_j$ in $D$ tending to $bD$ in
the sense that the $C_j$ eventually surround each compact subdomain of $D$,
such that 
$$
\int_{C_n}|g(z)|^p|dz|\leq M<\infty
$$ 
\endproclaim
The norm of $g\in E^p(D)$ is defined as 
$$
||g||^p_{E^p(D)}=\inf \sup_j\int_{C_j}|g(z)|^p|dz|
$$ 
where the inf is taken over all sequences of rectifiable curves $C_j$ in $D$ 
tending to $\partial D$.
\proclaim {Lemma 3.2} The holomorphic function  $F$ is in the Hardy
space $E^p(U)$.\endproclaim
\demo {Proof} Define a function $h$ on $\partial U$ by setting it
as  $$h(x+i\varphi(x,b))=f(x,b),\qquad h(x+i\varphi(x,0))=bf(x)$$ Observe that
$F\in E^1(U)$ if there is $H\in E^1(U)$ such that almost everywhere on
$\partial U$, the nontangential limit of $H$ equals $h$. Indeed in that case,
by Privalov's theorem, $H$ will agree with $F$.  According to Theorem 10.4 in
[Du], a necessary and sufficient condition for the existence of such $H$ 
is that $$\int_{bU}z^nh(z)dz=0,\quad n=0,1,2,\dots \tag 3.3$$ In our case, it
is clear that (3.3) will hold if we show it holds for $n=0$. The case $n=0$ is
equivalent to showing that $$\int_{bA}f(x,y)dZ(x,y)=0 \tag 3.4$$ where
$Z(x,y)$ is the first integral of $L$ defined above, and
$$A=(\alpha,\beta)\times (0,b)$$ In the integral (3.4), we have set
$f(x,0)=bf$. Note next that (3.4) holds if in the path $bA$ the segment $
(\alpha ,\beta)\times 0$ is replaced by $ (\alpha,\beta)\times \{\epsilon \}$,
for some $\epsilon >0$. In particular, this means that the integrals
$$\int_{(\alpha,\beta)}f(x,\epsilon)dZ$$ have a limit as
$\epsilon \to 0$.  Moreover, we also conclude that (3.4) will hold provided
that $$\lim_{\epsilon\to 0}
\int_{(\alpha,\beta)}f(x,\epsilon)dZ=\int_{(\alpha,\beta)}bf(x)\,dZ\tag 3.5$$
Thus in order to prove that $F\in E^1(U)$, it is sufficient to prove (3.5).
Since $\varphi(\alpha,t)$ is constant for $t>0$, we know that
$f(\alpha,t)$ is constant. Without loss of generality, we may assume that
$f(\alpha,t)\equiv 0$. Define then $g(x,t)=f(x,t)$ when $x>\alpha$ and
$g(x,t)=0$ when $x\leq \alpha$. The function $g$ is continuous, $Lg=0$ and for
$K$ compact in $x-$ space, $$\int_K|g(x,t)|\,dx=O(t^{-N})\quad \text {as}\quad
t\to 0^+$$ By Theorem 1.1, $\lim_{t\to 0^+}g(x,t)=bg$ exists in the
sense of distributions on $(-a,a)$. We will next show that
$$bg=\chi_{\alpha}(x)bf(x)$$ where $\chi_{\alpha}(x)$ denotes the
characteristic function of $(\alpha,a)$. Fix $\phi \in C_0^{\infty}(-a,a)$.
Recall from the proof of Theorem 1.1 that there is a smooth function
$\Phi ^N(x,t)$ such that $$\Phi ^N(x,0)=\phi (x),\quad |L^t\Phi ^N(x,t)|\leq
Ct^N$$ and that $$\langle bg,\phi \rangle = \int_{\alpha}^a
f(x,b)\Phi^{N}(x,b)\,dx-\int_0^b\int_{\alpha}^a f(x,t)L^t\Phi^N(x,t)\,dxdt \tag
3.6$$ For each $\epsilon >0$, let $\psi_{\epsilon}(x)\in C_0^{\infty}(\alpha
-\epsilon,a)$ such that $$(1)\quad \psi _{\epsilon}(x)\equiv 1 \text{ on supp
}\phi \cap (\alpha,a),\quad (2)\quad |D^k\psi_{\epsilon}(x)|\leq
\frac{c_k}{\epsilon^k}$$ Clearly, $$\int_{\alpha}^a
bf(x)\phi(x)\,dx=\lim_{\epsilon \to 0}\int_{-a}^a
bf(x)\psi_{\epsilon}(x)\phi(x)\,dx \tag 3.7$$ For the integrals on the right of
(3.7) we can use Theorem 1.1 to get: $$\align \int_{-a}^a
bf(x)\psi_{\epsilon}(x)\phi(x)\,dx &=\int_{-a}^a
f(x,b)\Phi_{\epsilon}^{N}(x,b)\,dx\\
&-\int_0^b\int_{-a}^af(x,t)L^t\Phi_{\epsilon}^N(x,t)\,dxdt\tag 3.8 \endalign $$
In the preceding expression, the function $\Phi_{\epsilon}^N(x,t)$ is chosen
using the proof of Lemma 3.3 in [BH1]. Indeed we recall from that Lemma that
for each $\epsilon >0$ and $k$ a nonnegative integer, there are smooth
functions $\phi^{\epsilon}_j(x,t)$ such that
$$\phi^{\epsilon}_0(x,t)=\psi_{\epsilon}(x)\phi(x)$$ and if
$$\Phi_{\epsilon}^k(x,t)=\sum_{j=0}^k\phi^{\epsilon}_j(x,t)\frac{t^j}{j!},$$
\roster 
\item $L^t\Phi_{\epsilon}^j(x,t)=\phi^{\epsilon}_{j+1}(x,t)\frac{t^j}{j!}$  and
\item $\phi^{\epsilon}_j(x,t)=\sum_{m\leq j} c_m(x,t)D_x^m\left(\psi_{\epsilon}(x)\phi(x)\right)$ 
\endroster 
where $c_m(x,t)$ are smooth functions independent of $\epsilon$ and
$c_m(x,t)=O((x-\alpha)^m)$.  We also have smooth functions $\phi_j(x,t)$ such
that $$\phi_0(x,t)=\phi(x)$$ and if
$$\Phi^k(x,t)=\sum_{j=0}^k\phi_j(x,t)\frac{t^j}{j!},$$
then
\roster 
\item $L^t\Phi ^j(x,t)=\phi_{j+1}(x,t)\frac{t^j}{j!}$  and
\item $\phi_j(x,t)=\sum_{m\leq j}
c_m(x,t)D_x^m\left(\phi(x)\right)$ 
\endroster

We consider now the first integral on the right in
(3.8).$$\align  \int_{-a}^a 
f(x,b)\Phi_{\epsilon}^{N}(x,b)\,dx&=\sum_{j=0}^N\int_{-a}^af(x,b) 
\phi_j^{\epsilon}(x,b)\frac{b^j}{j!}\,dx\\  &=\sum_{j=0}^N\sum_{m\leq 
j}\int_{-a}^af(x,b)c_m(x,b)D_x^m\left(\psi_{\epsilon}(x)\phi(x)\right) 
\frac{b^j}{j!}\,dx \endalign$$
$$=\sum_{j=0}^N\sum_{m\leq j}\sum_{k\leq  m}\binom
mk\int_{-a}^af(x,b)c_m(x,b)D_x^k\psi_{\epsilon}(x) 
D_x^{m-k}\phi(x)\frac{b^j}{j!}\,dx $$  In the terms above, when
$k>0$,  $$|c_m(x,b)D_x^k\psi_{\epsilon}(x)D_x^{m-k}\phi(x)|\leq C\epsilon 
^{m-k}$$  and the support of the integrand is contained in the interval 
$(\alpha -\epsilon, \alpha)$. Hence, as $\epsilon \to 0^+$, such terms go to 
$0$ while the term with $k=0$ converges to
$$\align
\sum_{j=0}^N\sum_{m\leq 
j}\int_{\alpha}^af(x,b)c_m(x,b)D_x^m\left(\phi(x)\right) \frac{b^j}{j!}\,dx&=
\sum_{j=0}^N\int_{\alpha}^af(x,b)  \phi_j (x,b)\frac{b^j}{j!}\,dx\\ &=
\int_{\alpha}^a  f(x,b)\Phi^{N}(x,b)\,dx \endalign $$

Therefore, we get $$\lim_{\epsilon \to
0} \int_{-a}^a f(x,b)\Phi_{\epsilon}^{N}(x,b)\,dx=\int_{\alpha}^a 
f(x,b)\Phi_{\epsilon}^{N}(x,b)\,dx\tag 3.9$$  which is the same as the first
term on the right in (3.6). We  next consider the second integral on the right
in (3.8).  $$
\int_0^b\int_{-a}^af(x,t)L^t\Phi_{\epsilon}^N(x,t)\,dxdt=\int_0^b\int_{-a}^a
f(x,t)\phi^{\epsilon}_{N+1}(x,t)\frac{t^N}{N!}\,dxdt $$
$$=\sum_{m \leq
N+1}\int_0^b\int_{-a}^a f(x,t)c_m(x,t)D_x^m\left(\psi_{\epsilon}(x) 
\phi(x)\right)\frac{t^N}{N!}\,dxdt$$
$$=\sum_{m\leq N+1}\sum_{k\leq
m}\int_0^b\int_{-a}^af(x,t)c_m(x,t)  \binom mk
D_x^k\psi_{\epsilon}D_x^{m-k}\phi\frac{t^N}{N!}  \,dxdt $$
 Again note that when $k>0$,
 $$|c_m(x,b)D_x^k\psi_{\epsilon}(x)D_x^{m-k}\phi (x)|\leq
 C{\epsilon}^{m-k},$$
 the $x-$ support of the integrand is contained in $(\alpha -\epsilon,
\alpha)$,  and since $f(.,t)=O(t^{-N})$, the term $f(x,t)t^N$ is bounded. It
 follows that as before, as $\epsilon \to 0^+$ such terms go to
 $0$ and we get:
$$\lim_{\epsilon \to
0}\int_0^b\int_{-a}^af(x,t)L^t\Phi^N(x,t)\,dxdt=
\int_0^b\int_{\alpha}^af(x,t)L^t\Phi_{\epsilon}^N(x,t)\,dxdt\tag 3.10$$
which is the same as the second integral on the right in (3.6).
>From (3.8), (3.9) and (3.10) we conclude:
$$\int_{\alpha}^abf(x)\phi(x)\,dx=\int_{\alpha}^af(x,b)\Phi^N(x,b)\,dxdt+
\int_0^b\int_{\alpha}^af(x,t)L^t\Phi^N(x,t)\,dxdt\tag 3.11 $$ We have thus
shown that $bg(x)=\chi_{\alpha}(x)bf(x)$ which implies that for any $\phi \in
C_c^{\infty}(-a,a)$,
$$\lim_{\epsilon \to
0}\int_{\alpha}^af(x,\epsilon)\phi(x)\,dx=\int_{\alpha}^abf(x)\phi(x)\,dx$$
Since $\varphi(\beta,t)\equiv 0$ for $t>0$, similar arguments imply that 
$$\lim_{\epsilon \to
0}\int_{\beta}^af(x,\epsilon)\phi(x)\,dx=\int_{\beta}^abf(x)\phi(x)\,dx$$
The preceding two limits establish (3.5) and hence (3.4). We have thus proved
that $F\in E^1(U)$. Hence by Theorem 10.4 in [Du], $F$ has a nontangential
limit $bF$ almost everywhere on $\partial U$ and that it can be
expressed as the Cauchy transform, $F=C^+(bF)$. But $bF=h\in L^p(\partial U)$.
Therefore, $F\in E^p( U)$. 
\enddemo

\heading {4. End of the proof of Theorem 3.1}
\endheading
\subhead {A. Bell-shaped regions and Hardy
spaces}\endsubhead
\bigskip

Consider a bounded region $\Omega\subset\ce$ satisfying the condition 
that there is $\alpha=\alpha (\Omega)>0$ with the property that almost every point $p$ in
the boundary admits a nonempty nontangential approach subregion
$$
\Gamma_{\alpha}(p)=\{z\in \Omega:|z-p|\leq
(1+\alpha)\text{dist}(z,\partial {\Omega})\}\tag 4.1
$$
that is, for a.e. $p \in\partial\Omega$, $\Gamma_\alpha(p)$ is open and 
$p$ is in the closure of $\Gamma_\alpha(p)$. This condition is 
satisfied, for instance, if  $\Omega$ is a bounded, simply connected 
region with rectifiable boundary. For this class of regions it is 
possible to define a class of  Hardy spaces as follows ([L]):

\proclaim {Definition 4.1}
Let $\Omega\subseteq \Bbb C$ be a bounded domain with a rectifiable 
boundary
and let $u$ be a function defined on $\Omega$. The nontangential
maximal function of $u$, $u^*$ and the nontangential limit of $u$,
$u^+$ are defined as follows:
$$
u^*(p)=\text{sup}_{\zeta\in\Gamma_{\alpha}(p)}|u(\zeta)|\quad
\text{a.e.}\quad p\in \partial \Omega,
$$
$$
u^+(p)=\lim_{\zeta\rightarrow p,\zeta\in
\Gamma_{\alpha}(p)}u(\zeta)\quad \text{a.e.}\quad p\in
\partial \Omega
$$
provided that the above limit exists.
\endproclaim

\proclaim{Definition 4.2}
For $1\le p<\infty$ the Hardy space is defined by
$$
H^p( \Omega)=\{f\in \Bbb O(\Omega):f^*\in L^p(\partial \Omega)\}
$$
where $\Bbb O(\Omega)$ denotes the holomorphic functions on $\Omega$.
\endproclaim

Our aim is to prove that $E^p(\Omega)=H^p(\Omega)$ for a particular class of domains $\Omega$ that includes the domain $U$ of Lemma 3.2. Let us point out that if $\Omega$ is the unit 
disc it is classical that both classes of Hardy spaces coincide and this fact implies ---by the Riemann mapping theorem--- that the same happens when $\Omega$ has smooth boundary. More generally, it is proved in [L] that $E^p(\Omega)=H^p(\Omega)$ also holds if $\Omega$ has a 
Lipschitz boundary and $1<p<\infty$.

We now consider smooth regions $U$ that are bounded by two smooth curves $\C_1$ and $\C_2$ that cross each 
other at two points $A$ and $B$ where they meet at angles $0\le\theta(A),
\theta(B)<\pi$. If $\theta(A), \theta(B)>0$ then $U$ has a Lipschitz boundary
and by the result mentioned before we know that $E^p(U)=H^p(U)$ for $p>1$. Our
methods will show that this equivalence still holds when the values
$\theta(A)=0$, $\theta(B)=0$ and  $p=1$ are allowed. By a conformal map
argument we may assume that \roster \item $A=0$ and $B=1$; \item the part 
$\C_1$ in the boundary of $U$ is given by $[0,1]\owns  t\mapsto t$ \item the
part  $\C_2$ in the boundary of $U$ is given by $[0,1]\owns  t\mapsto
x(t)+iy(t)$ where $x(t), y(t)$ are smooth real functions such that 
$x(0)=y(0)=y(1)=0$, $x(1)=1$. \endroster
We first prove that $H^p(U)\subset E^p(U)$. We construct for large 
integral $j$ a curve $C_j$ as follows. To every point 
$z\in\C_2\cap\partial U$ we assign the point 
$\gamma_{j,2}(z)=z+j^{-1}\bold n(z)$ where $\bold n(z)$ is the inward 
unit normal to $\C_2$   at $z$. For large $j$, $\C_2\owns 
z\mapsto\gamma_{j,2}(z)$ is a diffeomorphism and
$$
\text{dist}(\gamma_{j,2}(z),\C_2)=|\gamma_{j,2}(z)-z|=\frac{1}{j}.\tag 4.2
$$
Observe that the set
$$
D_j=\left\{z:\quad\text{dist}(z,[0,1]\times\{0\})\le\frac{1}{j} 
\right\}
$$
has a $C^1$ boundary $\partial D_j$ formed by 2 straight segments and 2 
circular arcs.
Fix a point $z_0\in\C_2$, choose  $j$ such that $z_0\notin D_j$ and 
consider the connected component of
$$
\left\{\gamma_{j,2}(z):\quad 
\text{dist}(\gamma_{j,2}(z),D_j)\ge\frac{1}{j}\right\}
$$
that contains $z_0$. Thus, we obtain a curve $C_{j,2}$ given by 
$[0,1]\supseteq[a_j,,b_j]\owns t\mapsto \gamma_{j,2}(x(t)+iy(t))\subset 
U$ that meets  $\partial D_j$ at its endpoints $A_j$, $B_j$ and remains 
off $D_j$ for $a_j<t<b_j$. Hence, we obtain a closed curve $C_j$ 
completing the curve $C_{j,2}$ with the portion $C_{j,1}$ of $\partial 
D_j$ contained in $U$ that joins $A_j$ to $B_j$. Because we are assuming 
that $\theta(A),\theta(B)<\pi$ we see that, for large $j$, $C_{j,1}$ is 
a horizontal segment at height $1/j$.
It is clear that all points in $C_j$ have distance $1/j$ to the 
boundary. Furthermore, if $q\in C_{j,2}$, $q\neq A_j$ and $q\neq B_j$ 
then $\text{dist}(q,\partial U)=\text{dist}(q,\C_2)=1/j$ because of (4.2) and
the fact that  $\text{dist}(q,[0,1]\times\{0\})>1/j$. Similarly, if $q\in
C_{j,1}$, $q\neq A_j$ and $q\neq B_j$ then $\text{dist}(q,\partial 
U)=\text{dist}(q,\C_1)=1/j$.
Thus, every point  $q\in C_j$ is at a distance $1/j$ of $\partial U$, we 
can always find $z\in\partial U$ such that 
$|q-z|=\text{dist}(q,\partial U)$  and $z$ is uniquely determined by 
$q$ except when $q=A_j$ or $q=B_j$ (in which case the distance may 
be attained at two distinct boundary points). In particular, whatever 
the value of $\alpha>0$, $q\in\Gamma_\alpha(z)$ for all $q\in C_j$ and 
$|g(q)|\le g^*(z)$ for any function $g$ defined on $U$.

Given $g\in H^p(U)$ we must show that
$$
\sup_j\int_{C_j}|g(z)|^p|dz|\le M<\infty.\tag 4.3
$$
We have
$$
\align
\int_{C_{j,2}}|g(q)|^p|dq|&=\int_{\gamma_{j,2}^{-1}(C_j)}
|g(\gamma_{j,2}(z))|^p\,|\gamma_{j,2}'(z)|\,|dz|\\
&\le \int_{\gamma_{j,2}^{-1}(C_j)}
|g^*(z)|^p\,|\gamma_{j,2}'(z)|\,|dz|\\
&\le C\int_{\C_2}|g^*(z)|^p\,|dz|.\tag 4.4
\endalign
$$
Similarly, using the map $\gamma_{j,1}(x)=x+i(1/j)\in C_{j,1}$, we get
$$
\int_{C_{j,1}}|g(q)|^p|dq|\le C\int_{\C_1}|g^*(z)|^p\,|dz|,\tag 4.5
$$
so adding (4.4) and (4.5) we obtain
$$
\int_{C_j}|g(q)|^p|dq|\le C \int_{\partial U}|g^*(z)|^p\,|dz|
$$
which implies (4.3) with $M=C\|g\|_{H^p}^p$.

To prove the other inclusion we first assume that $p=2$. Given $f\in 
E^2(U)\subset E^1(U)$ it has an a.e defined boundary value $f^+=bf\in 
L^2(\partial U)$ and the Cauchy representation
$$
f(z)=\frac{1}{2\pi i}\int_{\partial 
U}\frac{bf(\zeta)}{\zeta-z}\,d\zeta,\quad z\in U,
$$
is valid ([Du,Thm10.4]). Furthermore, $\|f\|_{E^p(U)}\simeq 
\|f^+\|_{L^p(\partial U)}$.

Following [L] we define
$$
T_*f^+(z)=\sup_{\epsilon>0}\left|\int_{\partial 
U\atop|\zeta-z|>\epsilon}
\frac{1}{\zeta-z}f^+(\zeta)\,d\zeta \right|,\quad z\in \partial U,
$$
and the maximal Hardy-Littlewood function
$$
Mf^+(z)=\sup\frac{1}{|I|}\int_I |f^+(\zeta)|\,|d\zeta|,\quad z\neq A,B
$$
where the $\sup$ is taken over all subarcs $I\subset\partial U$ that 
contain $z$ and $|I|$ denotes the arclength of $I$. Next we recall Lemma 
2.9 in [L]  that gives the estimate
$$
f^*(z)\le T_*f^+(z)+CMf^+(z),\quad z\in\partial U\setminus\{A,B\}.\tag 
4.6
$$
It is well known that $M$ is bounded in $L^2(\partial U)$. Furthermore, 
$T_*$ is also bounded in $L^2(\partial U)$ (this is a deep theorem when 
$U$ has just Lipschitz boundary ([C],[CMM])
but is much simpler here because $\partial U$ is smooth except at a 
couple of points). Then (4.6) implies that
$$
\|f\|_{H^2(U)}=\|f^*\|_{L^2(\partial U)}\le C\|f^+\| 
_{L^2(\partial U)}
\le C'\|f\|_{E^2(U)}.
$$
The same technique leads to the inclusion $E^p(U)\subset H^p(U)$ for 
$p>1$ because $T_*$ and $M$ are bounded as well in $L^p(\partial U)$ for 
$1<p<\infty$ but the method breaks down for $p=1$. So we recall that 
if $f\in E^p(U)$, $1\le p<\infty$, $f$ has a canonical factorization $f=FB$ where $F$ has no zeros, and $|B|\le1$. This is classical for 
the unit disc $\Delta$, where $B$ is obtained as a Blaschke product and 
the general case is obtained from the classical result. Indeed, if 
$w:\Delta\to U$ is a conformal map, it follows that $\tilde 
f(z)=f(w(z))(w'(z))^{1/p}$ is in $H^p(\Delta)$. Denote by $\tilde 
B(z)$ the Blaschke product associated to the zeros of $\tilde f$ counted 
with multiplicity. Then, $|\tilde B(z)|\le1$ has the same zeros as $f_1=f\circ w$ 
with the same multiplicity and if $0<r_j\nearrow 1$ it follows that
$$
\sup_j\int_0^{2\pi}\frac{|f_1(r_je^{i\theta})|^p}{|\tilde B(r_je^{i\theta})|^p}
|w'(r_je^{i\theta})|\,d\theta=\sup_j\int_0^{2\pi}| f_1(r_je^{i\theta})|^p
|w'(r_je^{i\theta})|\,d\theta\le C.\tag 4.7
$$

The proof of (4.7) is classical. It is clear that the supremum on the right of
(4.7) is bounded by the left hand side $\sup$, because $|\tilde B|\le1$. To
prove the reverse inequality one considers the finite product $\tilde B_N$ of
the first $N$ Blaschke factors. These partial products $\tilde B_N(z)\to\tilde 
B(z)$ normally in $\Delta$ as $N\to\infty$,  $|\tilde B_N(z)|=1$ for $|z|=1$ 
and $\tilde B_N$ is continuous on $|z|\le1$, so 
$$ 
\sup_j\int_0^{2\pi}\frac{|f_1(r_je^{i\theta})|^p}{|\tilde B_N(r_je^{i\theta})|^p}
|w'(r_je^{i\theta})|\,d\theta=\sup_j\int_0^{2\pi}| f_1(r_je^{i\theta})|^p
|w'(r_je^{i\theta})|\,d\theta 
$$ 
Then, using Fatou's lemma, 
$$ 
\align 
\sup_j\int_0^{2\pi}\frac{|f_1(r_je^{i\theta})|^p}{|\tilde B(r_je^{i\theta})|^p}
|w'(r_je^{i\theta})|\,d\theta&= \sup_j\int_0^{2\pi} \lim_{N\to\infty}
\frac{| f_1(r_je^{i\theta})|^p}{|\tilde B_N(r_je^{i\theta})|^p}
|w'(r_je^{i\theta})|\,d\theta\\ 
&\le\sup_j \liminf_{N\to\infty}\int_0^{2\pi} 
\frac{|f_1(r_je^{i\theta})|^p}{|\tilde B_N(r_je^{i\theta})|^p}
|w'(r_je^{i\theta})|\,d\theta\\ &\le\liminf_{N\to\infty}\sup_j \int_0^{2\pi} 
\frac{|f_1(r_je^{i\theta})|^p}{|\tilde B_N(r_je^{i\theta})|^p}
|w'(r_je^{i\theta})|\,d\theta\\
&\le\sup_j\int_0^{2\pi}|f_1(r_je^{i\theta})|^p|w'(r_je^{i\theta})|\,d\theta
\endalign
$$
Thus, if we set $B(\zeta)=\tilde B(w^{-1}(\zeta))$ we see that $|B|\le1$ in $U$, $F\doteq f/B$ does  not vanish in $U$, 
$F^{p/2}$ is well-defined  and (4.7) implies that $F\in E^p(U)$ (use as a
sequence of curves tending to $\partial U$ the images by $w$ of the circles of
radius $r_j$). Hence, $F^{p/2}\in E^2(U)$ and, by the case already proved,
$F^{p/2}\in H^2(U)$. This implies that $[F^{p/2}]^*=[F^*]^{p/2}\in
L^2(\partial U)$ so $F^*\in L^p(\partial U)$ and also, because $|B|\le1$,
$f^*=(FB)^*\in L^p(\partial U)$ which is what we wanted to prove. Summing up,
$E^p(U)=H^p(U)$ for $1\le p<\infty$.

\remark{Remark} We notice for later reference that the arguments above show  that any $f\in E^1(U)$ may be written as $f=g^2B$ with $|B|\le1$, $g\in E^2(U)$, $|bf|=|bg|^2$ and $\|f\|_{E^1}=\|g\|_{E^2}^2$.

\endremark

\subhead
B. Uniform bounds for traces
\endsubhead

We return to our solvable vector field with smooth coefficients 
 $$
L=\frac{\partial}{\partial y}+a(x,y)\frac{\partial}{\partial x}
$$
 We will assume without loss of generality that there is a smooth real function with compact support $\varphi(x,y)$ defined on $\erre^2$ such that $Z(x,y)=x+i\varphi(x,y)$ is a first integral of $L$, i.e., $LZ=0$ or, equivalently, $a(x,y)=-i\varphi_y(x,y)/(1+i \varphi_x(x,y))$. Furthermore, we will also assume that $\varphi(0,0)=\varphi_x(0,0)=0$ and that $|\varphi_x(x,y)|$ is
uniformly small throughout (this requirement will eventually become more  precise). Because $L$ satisfies the Nirenberg-Treves condition $(\P)$ ([NT],[T1]) it follows that
$$
\text{\it for every }x\in\erre \text{\it\ the map } \erre\owns y\mapsto\varphi(x,y) \text{\it\ is monotone.}
$$
We will also consider a homogeneous weak solution   $f(x,y)$ of class $C^0$ defined on  a rectangle
$$
(-a',a')\times (0,b')
$$ 
for some $a'>0$, $b'>1$ where it satisfies the equation $Lf=0$; we will also
assume that for $K$ compact, $\int_K f(x,y)dx$ has  tempered growth as
$y\searrow0$ and this implies that $f(x,y)$ possesses  a weak boundary value
at $y=0$ that will be denoted by $bf(x)$. We assume that the boundary value
$bf$ is in $L^p(-a',a')$ for some $1\le p\leq \infty$.  We wish to explore
whether the norms of the traces $f(.,y)$ in $L^p[-a,a]$ are bounded uniformly
in $y$ for some $0<a<a'$. The cases $p=1,\infty$ will be handled separately so
we assume henceforth that $1<p<\infty$.

Consider the graphs $\C_0$, $\C_1$ of the functions $y=\varphi(x,0)$ and
$y=\varphi(x,1)$ respectively. They cross at the origin and may cross many
more times in the strip $|x|<a'$. If they cross again for $|x|<a'$ to the
right and to the left of $x=0$ we restrict our attention to the interval
$[-a_1, a_2]\subset[-a',a']$ where $0<a_1, a_2<a'$ satisfy $\varphi(-a_1,0)=
\varphi(-a_1,1)$, $\varphi(a_2,0)= \varphi(a_2,1)$. We will assume initially 
that $a_1$ and $a_2$ exist. Define  $$ F=\left\{  x\in[-a_1,a_2]:\quad
\varphi(x,0)=\varphi(x,1)\right\}, $$ Recall from Lemma 2.3 that $(0,1]\owns
y\mapsto f(x,y)$ is independent of $y$ for a.e. $x\in F$. Thus, $$ \int_F
|f(x,y)|^p\,dx=\int_F |f(x,1)|^p\,dx \tag 4.8 $$ and we only need to study the
integrals $$ \int_{[-a_1,a_2]\setminus F} |f(x,y)|^p\,dx.
$$
Since $F$ is closed, $[-a_1,a_2]\setminus F=\bigcup_k(\alpha_k,\beta_k)$ with $\varphi(\alpha_k,0)=\varphi(\alpha_k,1)$, $\varphi(\beta_k,0)=\varphi(\beta_k,1)$ and
$\varphi(x,0)\neq\varphi(x,1)$ for $\alpha_k<x<\beta_k$. 
Thus, the curves $\C_0$ and $\C_1$ cross at the points
$A_k=(\alpha_k,\varphi(\alpha_k,0))$, $B_k=(\beta_k,\varphi(\beta_k,0))$ and
determine a region $U_k\subset\{x+iy:\alpha_k<x<\beta_k,
\varphi(x,0)<y<\varphi(x,1)\}$ between them. From Lemma 3.2 and the fact that
$|\varphi_x|$ is sufficiently small there is a function $F_k\in E^p(U_k)$ such
that $f=F_k\circ Z$ for $\alpha_k<x<\beta_k$, $0<y<1$. Furthermore the
boundary values of $F_k$ are given by $bF_k(x+i\varphi(x,0))=bf(x)$,
$bF_k(x+i\varphi(x,1))=f(x,1)$.

We  forget momentarily the region $U_k$ and consider the Cauchy transforms associated to $\C_0$ and $C_1$:
$$
C_j u(x)=\frac{1}{2\pi i}\int_{-\infty}^\infty\frac{1+i\varphi_x(x,j)}{x-x'+ 
i(\varphi(x,j)-\varphi(x',j))}u(x')\,dx',\quad j=0,1. \tag 4.9 $$
 
It is a celebrated  and deep theorem, first proved in full generality in [CMcM], that $C_j$ is bounded in $L^2$ if $x\mapsto\varphi(x,j)$ is Lipschitz, although we should mention that today  there exist rather short and elementary proofs of it ({\it cf.} [CJS]) and anyway the result is fairly simple when $x\mapsto\varphi(x,j)$ has bounded derivatives up to order 2 (in our case all derivatives are bounded) as it can be reduced to the continuity of the standard Hilbert tranform. It is a general fact concer
ing operators  associated to standard kernels ---as is the case of $C_j$--- that the operator norm of the truncated operators
$$
C_{j,\epsilon} u(x)=\frac{1}{2\pi i}\int_{|x-x'|>\epsilon}\frac{1+i\varphi_x(x,j)}{x-x'+ 
i(\varphi(x,j)-\varphi(x',j))}u(x')\,dx',\quad j=0,1. 
$$
is uniformly bounded in $\L(L^2(\erre))$ for  $0<\epsilon<\infty$. The
truncated operators are also associated to standard kernels (with bounds
uniform in $\epsilon>0$) which implies their uniform $L^p$-continuity,
$1<p<\infty$. Furthermore, the maximal operator $$ T_{j,*}u(x)=\frac{1}{2\pi
}\sup_{\epsilon>0}\left|\int_{|x-x'|>\epsilon}\frac{1+i\varphi_x(x,j)} {x-x'+
i(\varphi(x,j)-\varphi(x',j))}u(x')\,dx'\right|,\quad j=0,1,\tag 4.10 $$ is
also bounded in $L^p$, $1<p<\infty$, by a classical result of Calderón, Cotlar
and Zygmund (see, e.g., [CM,ch.IV]). We may also consider the Hardy-Littlewood
maximal operators associated to $\C_0$, $\C_1$, $$
M_ju(x)=\sup\frac{1}{\ell_j(I)}\int_I|u(x')|\,\sqrt{1+\varphi_x(x',j)^2}\,dx',\quad j=0,1 $$
where the $\sup$ is taken over all intervals that contain $x$ and 
$$
\ell_j(I)=\int_I \sqrt{1+\varphi_x(x,j)^2}\,dx,\quad j=0,1.
$$ 
Fixing $j$ and $p=x+\varphi(x,j)\in \C_j$ we consider regions of approach above and below $\C_j$ 
$$
\Gamma_{j,\alpha}(p)=\{z\notin \C_j:|z-p|\leq
(1+\alpha)\text{dist}(z,\partial {\C_j})\}\tag 4.11
$$ 
and by the arguments that led to (4.6), we have an analogous estimate 
$$
\sup_{z\in \Gamma_{j,\alpha}(p)}\left| \frac{1}{2\pi i}\int_{\C_j}\frac{u(\zeta)}{\zeta-z}\,d\zeta\right|\le T_{j,*}(u\circ\Pi_j^{-1})(x)+ C(\alpha)M_j(u\circ\Pi_j^{-1})(x)
\tag 4.12
$$
where $u \in L^p(\C_j)$, $C(\alpha)$ depends only on the aperture $\alpha$ and
$\Pi_j^{-1}$ is the inverse of the projection $\Pi_j:\C_j\to\erre$ given by
$x+i\varphi(x,j)\mapsto x$. We fix $\alpha$ from now on so that for
all $p\in \C_j$, $j=0,1$, the vertical line passing trough $p$ is contained in
$\Gamma_{j,\alpha}(p)$ with the obvious exception of the point $p$ itself.
This is possible because the curves that  bound $\Gamma_{j,\alpha}(p)$ meet at
$p$ forming with the normal at $p$ an angle $\theta=\pm
\cos^{-1}(1/(1+\alpha))$, and the normal is bounded away from the horizontal
direction.

Returning to the region $U_k$ and the function $F_k\in E^p(U_k)$, we may
represent $F_k$ as a Cauchy integral in terms of its boundary values $bF_k$ 
on $\partial U_k$. We may write  $bF_k=bF_{0,k}+bF_{1,k}$ according to the
decomposition $\partial U_k=(\partial U_k\cap\C_0) \cup (\partial
U_k\cap\C_1)$. We may extend $bF_{j,k}$ to $\C_j$, $j=0,1$, setting it equal
to zero off $ \partial U_k\cap\C_j$ so it becomes an element of $L^p(\C_j)$
with compact support.

Consider a point $q\in U_k$ and let $q_0$ and $q_1$ be the points in $\C_0$ and
$\C_1$ respectively that lie above $q$, i.e., $\Re q=\Re q_0=\Re q_1$. Our
choice of $\alpha$ then gives $$ q\in \Gamma_{0,\alpha}(q_0)\cap
\Gamma_{1,\alpha}(q_1). $$
Assuming that $\C_1$ is above $\C_0$ on $U_k$ we have
$$
2\pi i F_k(q)=\int_{\partial U_k}\frac{bF_k(\zeta)}{\zeta-q}\,d\zeta=
\int_{\C_0}\frac{bF_{0,k}(\zeta)}{\zeta-q}\,d\zeta-
\int_{\C_1}\frac{bF_{1,k}(\zeta)}{\zeta-q}\,d\zeta.\tag 4.13
$$
We may invoke (4.12) with $j=0$ and $p=q_0$ to estimate the first integral on
the right hand side of (4.13) by $2\pi$ times $$
T_{0,*}(b F_{0,k}\circ\Pi_0^{-1})(\Re q_0)+ C(\alpha)M_0(b F_{0,k}\circ\Pi_0^{-1})(\Re q_0)
\tag 4.14
$$
and there is an analogous estimate for the second integral. Fix now $0<y<1$, $\alpha_k<x<\beta_k$ and 
take $q=x+i\varphi(x,y)$ so (4.14) reads
$$
T_{0,*}[bF_k(x+i\varphi(x,0))](x)+ C(\alpha)M_0[bF_k(x+i\varphi(x,0))](x).
$$
Summing up, for $\alpha_k<x<\beta_k$, we have an estimate 
$$
|F_k(x+i\varphi(x,y))|\le K_0[bF_k(x+i\varphi(x,0))](x)+
K_1[bF_k(x+i\varphi(x,1))](x), $$
where $K_0$ and $K_1$ are bounded operators in $L^p$  independent of $k$.
Thus,
$$
\align
\int_{[-a_1,a_2]\setminus F} |f(x,y)|^p\,dx&=\sum_k \int_{\alpha_k}^{\beta_k} |f(x,y)|^p\,dx
=\sum_k\int_{\alpha_k}^{\beta_k} |F_k(x+i\varphi(x,y))|^p \,dx\\
&\le
C\sum_k\int_{\alpha_k}^{\beta_k}\left
(|K_0[bf](x)|^p+|K_1[f(x,1)](x)|^p\,dx\right )\\ &\le C_p\sum_k
\int_{\alpha_k}^{\beta_k}\left (|bf(x)|^p+|f(x,1)|^p\,dx\right ).\tag 4.15
\endalign $$ Thus, (4.8) and (4.15) yield
$$
\int_{-a_1}^{a_2} |f(x,y)|^p\,dx\le C\,\int_{-a_1}^{a_2}
\left (|bf(x)|^p+|f(x,1)|^p\,dx\right ),\quad 0<y<1. \tag 4.16
$$
So far we have assumed that the curves $\C_0$ and $\C_1$ cross each other both to the right and to the left of the origin. 
To finish the proof of the case $1<p<\infty$ we observe that if, say, $\C_0$
and $\C_1$ never cross to the right of the origin and we obtain a
region $U_1$ given by $0<x<a'$, $\varphi(x,0)<y<\varphi(x,1)$, we may
consider instead the  subregion $\tilde U_1$ given by $0<x<a''$,
$\varphi(x,0)<y<\varphi(x,1)$, for some $0<a''<a'$. An application of the
pointwise convergence result Theorem 2.1' to  $L$ on $x>0$ shows that for a.e.
$a''$, the function $(0,1)\owns y \mapsto f(a'',y)$ is bounded and hence in
$L^p(0,1)$. Finally, we may modify further $\tilde U_1$ by smoothing out the
corner in a neighborhood of the point $(a'',\varphi(a'',1))$ obtaining a
bounded region with 2 cusps of the type considered before and carry out our
analysis there. The details are left to the reader.

\subhead
C. The cases $\boldkey p$=1 and $\boldkey p$=$\boldsymbol\infty$
\endsubhead

Assume $p=1$. Keeping the notation of the previous subsection, it is clear that (4.8) holds for $p=1$ and we need only estimate 
$$
\int_{[-a_1,a_2]\setminus F} |f(x,y)|\,dx=\sum_k \int_{\alpha_k}^{\beta_k} |f(x,y)|\,dx
$$
uniformly in $y$. Using once again Lemma 3.2 we find functions $F_k\in
E^1(U_k)$ such that $f=F_k\circ Z$ on $(\alpha_k,\beta_k)\times(0,1)$ and
$bF_k(x+i\varphi(x,0))=bf(x)$, $bF_k(x+i\varphi(x,1))=f(x,1)$. In view of the
remark at the end of subsection A there are holomorphic functions $G_k\in
E^2(U_k)$ and  $B_k$ bounded by 1 in $U_k$ such that $F_k=G^2_k B_k$ and
$|bF_k|=|bG_k|^2$ a.e. The boundary of $U_k$ is bounded by the graphs of
$\varphi(x,0)$ and $\varphi(x,1)$ so let us denote by $b_0$ and $b_1$ the
corresponding boundary operators. Hence, invoking once more the operators
$K_0$ and $K_1$ which are continuous in $L^3$, we get $$ \align
\int_{\alpha_k}^{\beta_k} |f(x,y)|\,dx&=\int_{\alpha_k}^{\beta_k}
|F_k(x+i\varphi(x,y))|\,dx \le\int_{\alpha_k}^{\beta_k}
|G_k(x+i\varphi(x,y))|^2\,dx\\ &\le
C\int_{\alpha_k}^{\beta_k}\left (|K_0[b_0G_k](x)|^2+
|K_1[b_1G_k](x)|^2\,dx\right )\\ &\le
C_2\int_{\alpha_k}^{\beta_k}\left (|b_0G_k(x+i\varphi(x,0))|^2+
|b_1G_k(x+i\varphi(x,1))|^2\,dx\right )\\
&=C_2\int_{\alpha_k}^{\beta_k}\left (|b_0F_k(x+i\varphi(x,0))|+
|b_1F_k(x+i\varphi(x,1))|\,dx\right )\\
&=C_2\int_{\alpha_k}^{\beta_k}|bf(x)|+|f(x,1)|\,dx. \endalign $$ Thus, (4.15)
is also valid for $p=1$ and so is (4.16) which takes care of the case
$p=1$.\newline Finally we discuss the case $p=\infty$. By Lemma 2.3 $y\mapsto
f(x,y)$ is independent of $y$ for a.e. $x\in F$,  and so by  Lemma 2.2,
$bf(x)=f(x,1)$. Since $x\mapsto f(x,1)$ is continuous it follows that $bf(x)$
is essentially bounded on $F$. The complement of $F$ is a union of intervals
$(\alpha_k,\beta_k)$ associated to regions $U_k$ as described before. Since
$L^\infty[-a,a]\subset L^1[-a,a]$, the case $p=1$ implies that $f=F_k \circ Z$
on $U_k$  for some $F_k\in E^1(U_k)$ and (assuming that $\varphi_y\ge0$ for
$\alpha_k<x<\beta_k$) $bF_k$ is respectively given by $bf(x+i\varphi(x,0))$
and $f(x+i\varphi(x,1))$ on the two graphs that bound $U_k$. Thus, the
boundary value $bF_k$ of $F_k$ is essentially bounded by
$M=\|bf\|_{L^\infty[-a,a]}+\|f(\cdot,1)\|_{L^\infty[-a,a]}$ which implies that
$F_k$ itself is bounded by $M$ by the generalized maximum principle. Thus
$|f(x,y)|\le M$ for $\alpha_k<x<\beta_k$. Since $k$ is  arbitrary we conclude
that $M$ is a bound for $f(x,y)$.
\medskip Combining Corollary 1.3 and Theorem 3.1 for $p=1$ we get the
following: 
\proclaim 
{Corollary 4.1}  
Suppose $L=\frac{\partial}{\partial
t}+a(x,t)\frac{\partial}{\partial x}$ is locally solvable in $Q=(-a,a)\times
(-b,b)$, $f(x,t)$ a continuous function in $(-a,a)\times (0,b)$ such that
$Lf=0$ and for some $N$, 
$$ 
\int_K|f(x,t)|dx=O(t^{-N}) 
$$ 
on compact subsets $K$ of $(-a,a)$. Then the following are equivalent: 
\roster 
\item"(a)" $bf(x)$ is a locally finite measure.  
\item"(b)"  $bf(x)$ is a locally integrable function.
\item"(c)"  $N$ can be taken equal  to $0$ for any $K\subset\subset(-a,a)$.
\endroster 
\endproclaim
\demo {proof} Indeed, by Corollary 1.3, (a) implies (b), and Theorem 3.1
tells us that (b) implies (c). Finally, if (c) holds, then we can apply
Banach-Alaoglu's theorem to deduce (a).
\enddemo

\heading {5. Convergence regions for locally solvable vector fields}
\endheading

Suppose $\Omega$ is a smooth planar domain, $L=X+iY$ a {\it locally solvable} vector
field defined near each point of a compact connected portion $\Sigma$ of the boundary to which it is transversal, $f$ continuous on $\Omega$, $Lf=0$ where $L$ is defined, and for some defining function  $\rho$, there exists an integer $N$ such that the line integrals  
$$ \int_{\rho=t}|f|d\sigma_t =\quad O(t^{-N}),$$
where $\sigma_t$ denotes arc length on the curve $\rho=t$.
By Theorem 1.1, and the fact that  $\Sigma$ is noncharacteristic, we
know that $f$ has a boundary trace $bf$. Assume that this boundary value $bf\in
L^p(\Sigma)$ for some $1\le p\le\infty$. Since $\Sigma$ is noncharacteristic for $L=X+iY$,
by multiplying by $i$ if necessary, we may assume that $X$ is transversal
to $\Sigma$ and points  toward $\Omega$. For each $q\in \Sigma$, consider the integral curve $\gamma_q$ of $X$ through $q$ and its positive half  
$\gamma_q^+$ which enters $\Omega$. We shall distinguish below between two
types of points  $q\in\Sigma$: \roster \item "(I)" There exists a positive arc
$\{\gamma_q^+(\tau):\,\,0<\tau<\epsilon\}$ along which $X$ and $Y$ are
linearly dependent. \item "(II)" There is a sequence of points $q_k\in
\gamma_q^+$  converging to $q$ such that the vector field $L$ is elliptic at
each point of the sequence. \endroster We wish to attach to every point
$p\in\Sigma$ a subset $\Gamma(p)\subset\Omega$ such that: \roster \item $p$ is
an accumulation point of $\Gamma(p)$; \item if $q\in\Sigma$ is an accumulation
point of $\Gamma(p)$ then $q=p$; \item $\Gamma(p)$ contains an arc of
$\gamma^+_p$; \item for a.e. $p\in\Sigma$ 
$$
\lim_{\Gamma(p)\owns q \to p}f(q)=bf(p).
$$ 
\endroster
Since we are only interested in the behavior of $\Gamma(p)$ in arbitrary small
neighborhoods of  $p$, it would be more appropriate to consider the germ of
$\Gamma(p)$ at $p$ as well as of other related sets like $\gamma_p^+$ or
$S_{a,b}(p)$ defined below. However, to simplify the notation, this will be
done only implicitly and we shall not distinguish between sets and their
germs. It will be enough to carry out the construction of $\Gamma(p)$ for $p$
in a small neighborhood in $\Sigma$ of a given point of $\Sigma$. In order to
define $\Gamma(p)$, fix $p\in \Sigma$ and consider a first integral $Z$ of $L$
defined in a neighborhood of $p$, such that the restriction of $\Re Z$ on
$\Sigma$ has a nonzero differential. For $(a,b)$ a real vector close to $(1,0
)$, define  $$ S_{a,b}(p)=\{w:\Re (a+ib)Z(w)=\Re (a+ib)Z(p)\} $$  For $(a,b)$
in a small disk $V\subset\erre^2$ centered at $(1,0)$ these curves are
transversal to $\Sigma$. Let $S_{a,b}^+(p)$ denote the part of $S_{a,b}(p)$
that is in $\Omega$ and set $$ \Gamma(p)=\bigcup_{(a,b)\in V} S_{a,b}^+(p). $$
We now discuss whether $\Gamma(p)$ enjoys properties (1) through (4). It is
clear by construction that (1) and (2) are satisfied because $\Gamma(p)$ is a
union of curves entering $\Omega$ transversally. To check (3) we observe that if $p$ is of 
type (I) $S_{a,b}^+(p)\subset\gamma^+_p$, $(a,b)\in V$, in particular,
$\Gamma(p)=S_{1,0}^+(p)=\gamma^+_p$ locally. If $p$ is of type (II),  we may
choose the coordinates so that in a neighborhood of the origin
$\Sigma=\{t=0\}$, $\Omega=\{t>0\}$, $Z(x,t)=x+i\phi(x,t)$,
$\phi(0,0)=\phi_x(0,0)=0$, $X=\partial_t-[\phi_x\phi_t/(1+
\phi_x^2)]\partial_t$ and $p=(x_0,0)$.  Hence, $\gamma^+_p$ can be
parametrized as  $(x_1(s),s)$ where $x_1$ satisfies the ODE $$
\frac{dx_1}{dt}=-\frac{\phi_t\phi_x}{1+\phi_x^2},\quad x_1(0)=x_0 $$ while
$S^+_{a,b}(p)$ is given by the graph of $x=x_2(t)$ where $x_2$ satisfies the
implicit equation $$ x_2=\frac{b}{a}\phi(x_2,t),\quad x_2(0)=x_0. $$ We now
look at the images of $S^+_{a,b}(p)$ and $\gamma^+_p$ under the map
$(x,t)\mapsto(\xi,\eta)$, $\xi=x$, $\eta=\phi(x,t)$, and call this images
$\tilde S^+_{a,b}(p)$ and $\tilde \gamma^+_p$ respectively. With a slight
abuse of notation this map can be denoted by $Z$. So $\tilde
S^+_{a,b}(p)=Z(S^+_{a,b}(p))$ is an interval of the line
$a(\xi-x_0)-b(\eta-\eta_0)=0$, $\eta_0=\phi(x_0,0)$, and
$\tilde\gamma^+_p=Z(\gamma^+_p)$ may be parametrized as $(\xi(\eta),\eta)$
where $\xi=\xi(\eta)$ is seen to satisfy, after a short computation: $$
\frac{d\xi}{d\eta}=-\phi_x;\quad \xi(\eta_0)=x_0.\tag 5.1 $$ Since $L$
satisfies condition $(\P)$ and the origin is of type (II), $\phi$ has a
consistent sign in a neighborhood of the origin, say $\phi\ge0$. From standard
estimates for positive functions ([Di], [Gl]) it follows that $|\phi _x|\le
C\sqrt{\phi}$. Thus, (5.1) shows that $\tilde\gamma_p^+$ satisfies the
differential inequality $$ \left|\frac{d\xi}{d\eta}\right|\le
C\sqrt{\eta};\quad \xi(\eta_0)=x_0. $$ Therefore, it is contained in the
sector bounded by the straight lines $\xi=x_0\pm b(\eta-\eta_0)$,
$\eta>\eta_0$, for any positive $b$ and  $\eta_0<\eta\le\eta_1(b)$ if
$\eta_0(x_0)$ and $\eta_1(b)$ are taken small enough, and it follows that 
$\tilde\gamma_p^+$ is contained in the union $\bigcup\tilde S^+_{a,b}(p)$, or 
$$ Z(\gamma_p^+)\subset \bigcup_{(a,b)\in V} Z\left(S^+_{a,b}(p)\right). $$ On
the other hand, since $t\mapsto\phi(x,t)$ is monotone, the inverse image of a
point $q=(\xi,\eta)$, $Z^{-1}\{q\}$, is of the form $\{\xi\} \times [c,d]$
where $\phi_t(\xi,\eta)=0$ for $\eta\in[c,d]$. It is then easy to conclude
that if $Z^{-1}\{q\}$ intersects $\gamma^+_p$ (resp. $S^+_{a,b}(p)$) it is
totally contained in $\gamma^+_p$ (resp. $S^+_{a,b}(p)$) which implies that
$\gamma^+_p=Z^{-1}(Z(\gamma^+_p))$ and $S^+_{a,b}(p)=Z^{-1}(Z(S^+_{a,b}(p)))$.
Then the above inclusion implies that a small arc of $\gamma_p$ is contained
in $\bigcup S^+_{a,b}(p)$. 

Next we discuss the validity of (4). First, we point out that if $p$ is of
type (I) then $\Gamma(p)=\gamma^+_p=S^+_{1,0}(p)$ in a neighborhood of $p$ and
we may apply Theorem 2.1 to obtain the desired convergence result. More
generally, 
\proclaim  
{Theorem 5.1} For almost all $p \in \Sigma$  $$
\lim_{\Gamma(p)\owns q \mapsto p}f(q)=bf(p).
$$
\endproclaim 
\demo{Proof} 
We may assume that $p$ is of type (II) by the preceding comments. The
hypotheses tell us that $x=\Re Z$ and $t=\rho$ form a change of coordinates
near a point, say $p\in \Sigma$. We may assume $p$ is mapped to the origin. In
these new coordinates, $\Sigma$ is mapped to $t=0$,  $L$ takes the form $$
\frac{\partial}{\partial t}+a(x,t)\frac{\partial}{\partial x}
$$
except for a nonvanishing factor. The first integral $Z(x,t)=x+i\phi(x,t)$. We
now recall that for some rectangle $Q_r=(-r,r)\times (0,r)$, there is a
holomorphic function $F\in H^p(U)$, $U=Z(Q_r)$, such that $f(x,t)=F(Z(x,t))$.
We focus on the boundary piece of $U$ given by $\Sigma_0=Z(x,0)$. We know that
there is $\alpha >0$ such that if 
$$\widehat{\Gamma}_{\alpha}(q)=\{z\in U:|z-q|\leq (1+\alpha)d(z,\partial U)\}$$
then for almost all $q\in \Sigma_0$, 
$$\lim _{\widehat{\Gamma}_{\alpha}(q)\owns z\mapsto
q}F(z)=bF(q)$$ Fix $q=x_0+i\phi(x_0,0)$ where this limit exists. For $(a,b)\in V$, consider the curve 
$$S_{a,b}(x_0)=\{(x,t):ax-b\phi(x,t)=ax_0-b\phi(x_0,0)\}$$
Let $S_{a,b}^+(x_0)$ be the part of $S_{a,b}(x_0)$ where $t>0$. Observe that
the theorem will be proved if we show that $Z(S_{a,b}^+(x_0))$ is contained in
$\widehat{\Gamma}_{\alpha}(q)$. Let $(y,s)\in S_{a,b}^+(x_0)$. Then
$b(\phi(y,s)-\phi(x_0,0))= a(y-x_0)$, and so 
$$|y+i\phi(y,s)-q|^2=\left(1+\frac{b^2}{a^2}\right)|\phi(y,s)-\phi(x_0,0)|^2
\tag 5.2$$
If the point $(y,s)$ is close enough to $(x_0,0)$, there is $x$ such that
$$\align
d(y+i\phi(y,s),\partial U)^2
&=|y-x|^2+|\phi(y,s)-\phi(x,0)|^2\\ &\geq |y-x|^2+(1-\epsilon
^2)|\phi(y,s)-\phi(x_0,0)|^2+\\ & \left(1-\frac{1}{\epsilon ^2}\right
)|\phi(x,0)-\phi(x_0,0)|^2
\\ &\geq |y-x|^2+(1-\epsilon^2)|\phi(y,s)-\phi(x_0,0)|^2-2||\D_x\phi|||x-x_0|^2
\endalign $$
for $\epsilon$ close to $0$. Here $||D_x\phi||$ denotes the sup norm which may
be taken to be as small as we wish from the outset. We now consider two cases.
Suppose first $|y-x_0|\leq |y-x|$. Then $|x-x_0|\leq 2|y-x|$, and so 
$$\align
d(y+i\phi(y,s),\partial U)^2&\geq
|y-x|^2+(1-\epsilon^2)|\phi(y,s)-\phi(x_0,0)|^2-2||D_x\phi|||x-x_0|^2\\
&\geq |1-8||D_x\phi|||y-x|^2+(1-\epsilon^2)|\phi(y,s)-\phi(x_0,0)|^2
\tag 5.3
\endalign $$

Comparing (5.2) with (5.3), if we take $\epsilon$ small enough and $(a,b)$ is
sufficiently close to the vector $(1,0)$, then $Z(y,s)\in \widehat{\Gamma}_q$. Suppose
next $|y-x_0|\geq |y-x|$. Then $|x-x_0|\leq 2|y-x_0|$. Hence
$$\align
 d(y+i\phi(y,s),\partial U)^2&\geq
(1-\epsilon^2)|\phi(y,s)-\phi(x_0,0)|^2-8||D_x\phi|||y-x_0|^2\\
&=\left (1-\epsilon^2-8||D_x\phi||^2\frac{b^2}{a^2}\right
)|\phi(y,s)-\phi(x_0,0)|^2
\endalign
$$
 which again shows that by choosing $(a,b)$ close enough to $(1,0)$,
we get $Z(y,s)\in \widehat{\Gamma}_q$. The assertion of the theorem has thus been proved.
\enddemo
Summing up, we may think of the regions of convergence $\Gamma(p)$ as ``cusps" stemming from $p$ and entering $\Omega$ that contain $\gamma^+_p$. If $p$ is of type (I), $\Gamma(p)$ reduces to $\gamma_p^+$ but when $p$ is of  type  (II), $\Gamma(p)$ contains a neighborhood of $\gamma^+_p$ in $\Omega$.

\example{Example 5.1}
Let 
$$
L=\frac{\partial}{\partial t}-2it\frac{\partial}{\partial x},\quad \Omega=\{(x,t):\,\, t>0\},
\quad Z=x+it^2.
$$
Here, for any $p=(x_0,0)$, we may take $\Gamma(p)$ as the cusp bounded by two parabolas
$$
\Gamma(p)=\{(x,t):\quad x_0-ct^2<x<x_0+ct^2,\,\,t>0\}.
$$
\endexample
\heading {6. A uniqueness result and an application}
\endheading

We keep the notation of Section 5 and consider a  smooth planar domain $\Omega$,  
a locally integrable vector field $L=X+iY$ defined near each point of a 
closed subinterval $\Sigma$ of the boundary to which it is transversal, a
function $f\in C^0(\Omega)$ satisfying $Lf=0$ where $L$ is defined and such
that for some defining function  $\rho$, there exists an integer $N$ such that
the line integrals   $$  \int_{\rho=t}|f|d\sigma_t =\quad O(t^{-N}), $$
where $\sigma_t$ denotes arc length on the curve $\rho=t$. We thus know that $f$ has a 
trace $bf$ defined in a neighborhood of $\Sigma$. Assume that this 
boundary value is a  finite measure $\mu$. If we  denote by $\sigma$ the arc
length measure on $\partial\Omega$, Corollary 4.1 shows that $\mu$ is
absolutely continuous with respect to $\sigma$. A moment's reflection about
the example $L=\partial_t$, $\Omega=\{t>0\}$, shows that, in general, the
converse is not true, i.e., $\sigma$ need not be absolutely continuous with
respect to $|\mu|$ (the total variation of $\mu$), even if $\mu$ is not
identically zero. On the other hand, this phenomenon is not possible at points
where the behavior of $L$ at the boundary of $\Omega$ is removed from that of
a real vector field. More precisely, consider a Borel set $E\subset\Sigma$
such that $|\mu|(E)=0$ and let $p\in E$ be a point of type (II). Then, in
local coordinates, we may assume that $p=(0,0)$, $\Sigma=\{t=0\}$,
$\Omega=\{t>0\}$, $Z(x,t)=x+i\phi(x,t)$ is a local first integral satisfying
$\phi(0,0)=\phi_x(0,0)=0$. Since $(0,0)$ is of type (II), there is a sequence
$t_j\mapsto 0$ such that $\phi (0,t_j)\neq 0$. Withoult loss of generality, we
may assume that $\phi (0,t_j)>0$. We can then apply Theorem 3.1 in [BH1] as
in the proof of Theorem 2.1' to get a holomorphic function $F$ of tempered
growth defined on 
$$
Q=\{Z(x,0)+iZ_x(x,0)v:x\in (-r,r),\quad 0<v<\delta \}
$$
such that for any $\psi \in C_c^{\infty}(-r,r)$
$$
\int bf(x)\psi(x)dx=\lim_{v\mapsto 0}\int F(Z(x,0)+iZ_x(x,0)v)\psi(x)dZ(x,0)
$$
where $bf=\mu$. Since $bf$ is a locally integrable function, as is well known, the holomorphic 
function $F$ converges nontangentially in the region $Q$ to $bf(x)$ a.e.on the
part $\{Z(x,0)\}$ of the boundary of $Q$. Then, by the Riesz uniqueness
theorem, either $F$ vanishes identically in $Q$ or the zero set of $bF$ is a
subset of $\partial U$ with  null linear measure. In other words, there is a
neighborhood $V$ of $p$ in $\erre^2$ such that either $f\equiv0$ in
$V\cap\Omega$ or $\sigma(V\cap E)=0$. We now denote by $\Sigma_1\subset\Sigma$
(resp. $\Sigma_2\subset\Sigma$) the set of  points of $\Sigma$ of type (I)
(resp. of type (II)) and assume that \roster \item "$(*)$" {\it for any
$p\in\Sigma_2$ and any neighborhood $V$ in $\erre^2$ of $p$, $f$ does not
vanish identically on $V\cap\Omega$;} \endroster Then,  we have shown that if
$(*)$ holds  $|\mu|(E)=0$ implies that $\sigma(E\cap\Sigma_2)=0$ or,
equivalently, that $E\subset\Sigma_1$ except for a $\sigma$-null set. This can
be restated by saying that on $\Sigma_2$, $\sigma$ and $|\mu|$ are mutually
absolutely continuous with respect to each other. In fact, the  argument
shows more. Let's  recall that $p\in E\subset\Sigma$ is called a
$\sigma$-density point of $E$ if $\sigma(V\cap E)>0$ for any neighborhood $V$
of $p$. We have \proclaim   {Theorem 6.1}  Let $L$, $\Omega$, $f\in
C^0(\Omega)$ and $\Sigma$  as above and assume that  \roster \item $Lf=0$ on
$W\cap\Omega$ for some open $W\supset\Sigma$;  \item $bf\in L^1(\Sigma)$;
\item $p\in\Sigma_2$ is a $\sigma$-density point of the set $E=\{bf(x)=0\}$.
\endroster Then, there is an open disc $\Delta=\Delta(p,r)$ such that $f$
vanishes identically on $\Delta\cap\Omega$.
\endproclaim 
Suppose next that $L, \Omega, f\in C^0(\Omega)$ and $\Sigma$ are as above
except that we no longer make the growth assumption on the line integrals of
$|f|$. In particular, we don't assume that $f$ has a trace on $\Sigma$. We
then get the following convergence result generalizing to locally integrable
vector fields a classical result for holomorphic functions : 
\proclaim {Corollary 6.2} Assume that $\Re f \geq 0$. Then for almost all
$p\in \Sigma_2$, $\lim_{\gamma_p^+\ni q\mapsto p}f(q)$ exists and is finite. If
$L$ is locally solvable, the limit can be taken in the set $\Gamma (p)$.
\endproclaim 
In the corollary, we are using notations introduced in section 5.
\demo {Proof} Let $F=\frac{1}{1+f}$. Observe that $LF=0$, and $F$
is bounded. Therefore, we can apply Theorem 2.1 or the results in section 5 to
deduce convergence for $F$. Since $f=\frac{1}{F}-1$, if $p\in \Sigma_2$ is a
point of convergence for $F$, it is also a point of convergence for $f$,
unless $bF(p)=0$. Since $F$ does not vanish identically, by Theorem 6.1, such
points $p$ form a set of measure zero.
\enddemo

\enddocument